\RequirePackage{fix-cm}
\documentclass[smallextended]{svjour3}       
\smartqed  
\usepackage{amssymb,amsmath}
\usepackage{bbm}
\usepackage{dsfont}
\usepackage[english]{babel}
\usepackage{graphicx}
\usepackage{mathrsfs}
\usepackage{amsfonts}
\usepackage{fancyhdr}
\usepackage{textcomp}
\usepackage[official]{eurosym}
\usepackage{algorithm}
\usepackage{algorithmic}
\usepackage{geometry}
\usepackage{tikz}

\usepackage{epstopdf}
\DeclareGraphicsExtensions{.eps}

\definecolor{red}{rgb}{1,0,0}

\newcommand\pr{\prime}

\renewcommand\pr{\mathrm{Pr}}
\newcommand\xxi{{\boldsymbol{\xi}}}

\newlength{\figwidth}
\newlength{\tabwidth}

\newcommand\NN{{\hbox{I\kern-.14em{N}}}}
\newcommand\RR{{\hbox{I\kern-.14em{R}}}}
\newcommand\ZZ{{\hbox{I\kern-.14em{Z}}}}

 \journalname{  }
\begin{document}

\title{Sampling methods for multistage robust convex optimization problems\thanks{The work has been supported under Bergamo University grant 2015-2016 and by funds of the CNR-JST Joint International Lab COOPS.}
}


\author{Francesca Maggioni  \and
				Marida Bertocchi \and
				 Fabrizio Dabbene \and
				Roberto Tempo 
}


\institute{Francesca Maggioni, Bergamo University \at
              Via dei Caniana n. 2 24127 Bergamo, Italy \\
              Tel.: +39-0352052649\\
              Fax: +39-0352052549 \\
              \email{francesca.maggioni@unibg.it}           
           \and
           Marida Bertocchi, Bergamo University \at
              Via dei Caniana n. 2 24127 Bergamo, Italy
			      \and
						Fabrizio Dabbene, CNR-IEIIT, Politecnico di Torino \at
              C.so Duca degli Abruzzi n. 24
              10129 Torino, Italy
						\and
           Roberto Tempo, CNR-IEIIT, Politecnico di Torino \at
              Corso Duca degli Abruzzi  n. 24
           10129 Torino, Italy
}

\date{Received: date / Accepted: date}

\maketitle

\begin{abstract}
In this paper, probabilistic guarantees for constraint sampling of multistage robust convex optimization problems are derived. The dynamic nature of these problems is tackled via the so-called scenario-with-certificates approach. This allows to avoid the conservative use of explicit parametrizations through decision rules, and provides a significant reduction of the sample complexity to satisfy a given level of reliability. An explicit bound on the probability of violation is also given. Numerical results dealing with a multistage inventory management problem
show the efficacy of the proposed approach.
\end{abstract}


\section{Introduction} 
{In many practical situations, the decision process is affected by uncertainty. In such cases, the so-called \textit{uncertainty set} where all realizations of the random parameters lie is considered, and then optimize an objective function protecting
against the worst possible uncertainty realization. This is the key philosophy behind the robust optimization modeling paradigm.
}
The original robust optimization models deal with \textit{static} problems, where all the decision variables have to be determined before the uncertain parameters are selected.
A vast literature focused on uncertainty structure to obtain computationally tractable problems is available, see for instance \cite{10.2307/30036559}  and   \cite{10.2307/168933} for polyhedral uncertainty sets and  \cite{BenTal19991} for ellipsoidal uncertainty sets, respectively.

However, this approach cannot directly handle problems that are multiperiod in nature, where a decision at any period should take into account data realizations in previous periods, and the decision maker needs to adjust his/her strategy on the information revealed over time. This means that some of the variables (non-adjustable variables) must be determined before the realization of the uncertain parameters, while the other  variables (adjustable variables)  have to be chosen after the uncertainty realization. For a recent overview of multiperiod robust optimization, we refer to  \cite{doi:10.1137/080734510}, \cite{Delage_2015} and \cite{Gabreletal2012}.
In order to describe such  a situation, and extend  robust optimization to a dynamic framework,  the concept of \textit{Adjustable Robust Counterpart} (ARC) has been first  introduced and analyzed in  \cite{Ben-Tal2003}. This approach opened up several new research directions, such as portfolio optimization    \cite{Ceria2006}, {\cite{PiNar:2005:RPO:2308917.2309431}},   \cite{Tutuncu2004},  inventory management \cite{Ben-Tal:2005:RFC:1246360.1246366},\cite{doi:10.1287/opre.1050.0238}, scheduling   \cite{Lin20041069},  \cite{Yamashita2007}, facility location \cite{POMS:POMS1194},  revenue management   \cite{doi:10.1287/msom.1080.0252} and energy generation   \cite{6476049}.  ARC is clearly less conservative than the static robust approach, but in most cases it turns out to be computationally intractable.
One of the most recent methods to cope with this difficulty is obtained  by approximating the adjustable decisions by \textit{decision rules}, i.e. linear combinations of given basis functions of the uncertainty. A particular case  is the \textit{Affinely Adjustable Robust Counterpart} (AARC)  \cite{Ben-Tal2003}, where the adjustable variables are  affine functions of the uncertainty.
The decision rule approximation often allows to obtain a formulation which is equivalent to a tractable optimization problem  (such as linear, quadratic and second-order conic   \cite{doi:10.1137/S1052623401392354},  or semidefinite   \cite{doi:10.1137/S1052623496305717}), transforming the original dynamic problem into a static robust optimization problem whose decision variables are the coefficients of the linear combination.

However, in many practical cases, also the static robust optimization problem ensuing from the decision rule approximation is still  numerically intractable. In these situations, one can recur to approximate solutions based on constraint sampling, which consists in taking into account only a finite set of constraints, chosen at random among the possible continuum of constraint instances of the uncertainty.
 The attractive feature of this method is to provide  explicit bounds on the measure  of the original constraints of the static problem that are possibly violated by the randomized solution.
The properties of the solutions provided by this approach, called   scenario approach have been studied  in \cite{Calafiore2004}, \cite{doi:10.1137/07069821X}, \cite{deFarias:2004:CSL:1024112.1024150}, where it has been shown that most of the constraints of the original static problem are satisfied provided the number of samples sufficiently large.
The constraint sampling method has been also extensively studied within the chance constraint approach through different directions by \cite{Erdogan:2006:ACC:1124589.1124592}, \cite{doi:10.1137/070702928}, \cite{doi:10.1137/050622328} and \cite{Pagnoncelli2009}.


In  \cite{4434596}, \cite{1428665}, \cite{Vayanos2012459}, multistage convex robust optimization problems are solved by combining general nonlinear decision rules and constraint sampling techniques. This means that the dynamic robust optimization problem is transformed into a static robust optimization problem through decision rules approximation and then solved via a scenario counterpart.
In practice, the novelty of \cite{Vayanos2012459} is to  introduce, besides polynomial decision rules, also trigonometric monomials and basis functions based on sigmoidal and Gaussian radial functions, thus allowing  more flexibility.
A rigorous convergence proof for the optimal value, based on the decision rule approximation and of the constraint randomization approach is also investigated.  Convergence is proved when both the  complexity  parameter (number of basis in the decision rule approximation) and    the number of samples tends to infinity.


In the context of randomized methods for uncertain optimization control problems,  the \textit{scenario with certificates} approach  has been proposed in \cite{Formentin2015RobustSA},  based on an original idea of \cite{Oishi2006}. This approach has been then extended and exploited for anti-windup augmentation problems {\cite{Formentin2015RobustSA}}.
The main idea of this approach is to distinguish between design variables (corresponding to non-adjustable variables) and certificates (corresponding to adjustable variables).

In this paper, we  consider randomized methods for   robust convex multistage optimization problems.
We treat the dynamic nature of the problem  via the scenario with certificates approach, thus avoiding the conservative use of parametrization through decision rules. This implies a significant reduction of the number of samples required to satisfy the level of reliability of the constraints.
In particular, we show that a multistage robust linear optimization problem $\textrm{RO}_H$,  is equivalent to a linear robust optimization problem with certificates $\textrm{RwC}_H$, and a bound on the probability of violation is   provided for the scenario  with certificates problem  $\textrm{SwC}_H^N$ based on $N$  instances (or scenarios) of the uncertain constraints and $H$ stages.  The analysis has been extended to the convex case.
Furthermore, upper and lower bounds obtained by relaxing the nonanticipativity constraints are also provided.

The rest of the paper is as follows.
Section \ref{problemformulation}  discusses the formulations of two-stage, multistage robust linear and convex programs  and provides a result on  the probability of violation of constraints. Bounds on the number of scenarios needed to obtain a user-prescribed  guarantee of violation  is given. Section \ref{boundsMultistageRO} provides a  chain of inequalities among lower bounds on the multistage robust optimization problem. Section \ref{numericalresults}  presents several numerical results dealing with
a multistage inventory management problem. The conclusions follow.

\section{Problem formulation}
\label{problemformulation}
\subsection{Notation}
\label{notation}
In this paper, the uncertainty is described by a  discrete random process  $\xi^t$, $t=1,\dots, H$,  defined on a probability space $(\Xi^t,\mathscr{A}^{t},\pr)$  with marginal support $\Xi^t\subseteq \mathbb{R}^{n_t}$ and given probability distribution $\pr$ on the $\sigma-$algebra $\mathscr{A}^t$ (with $\mathscr{A}^t\subseteq \mathscr{A}^{t+1}$). The process $\xi^t$   is revealed gradually over discrete times in $H$ periods, and $\xxi^{t}:=(\xi^1,\dots,\xi^t)$, $t=1,\dots,H-1$  denotes the history of the process up to time $t$.

The decision variable at each discrete time  is indicated with  $x^t\in\mathbb{R}^{n_t}$, $t=1,\dots,H$.  The decision $x^1$ is selected at time (stage) $1$ before the future outcome of $\xi^1$ is revealed,
the decision  $x^{t}$ at stage  $t=2,\dots,H$ is  $\mathscr{A}^{t-1}$-measurable and it depends on the information up to time $t$. More precisely  the decision process is {\it nonanticipative}, i.e. it  has the form
\begin{eqnarray}
 &\textrm{decision}(x^1)&\rightarrow \textrm{observation}(\xi^1)\rightarrow \textrm{decision}(x^2)\rightarrow \textrm{observation}(\xi^2)\rightarrow \dots\nonumber\\  
 &\dots & \rightarrow \textrm{decision}(x^{t-1})\rightarrow \textrm{observation}(\xi^{t-1}) \rightarrow \textrm{decision}(x^{t})\rightarrow\dots\nonumber\\
&\dots & \rightarrow \textrm{observation}(\xi^{H-1}) \rightarrow \textrm{decision}(x^{H}).\nonumber
 \end{eqnarray}

In the following {\sf X} denotes the Cartesian product among sets, and the Binomial distribution with parameters $\epsilon\in\mathbb{R}$, $N,n\in\mathbb{N}$, $N>n$, is denoted as
$B(N,\epsilon,n+1)$.

\subsection{Two-stage robust linear case}
\label{sec:TwoStageRobustLinearProgram}
To simplify our exposition, we first analyze a simple  \textit{two-stage robust linear program}, formally defined as follows\footnote{We adopt the convention of putting as pedices the number of stages of the problem, e.g.  RO$_2$ denotes a two-stage robust linear problem.}
\begin{eqnarray}
\label{RO2}
\textrm{RO}_2 &:=&\! \min_{x^1} c^{{1}^{\top}}{x^1} + \sup_{\xi^{1}\in\Xi^1}\!\!\left[\min_{{x^2}(\xi^{1})} c^{{2}^{\top}}\left(\xi^{1}\right) {x^2}(\xi^{1})\!\right]  \\
& & \textrm{s.t. }			 A{x^1}=h^1,\quad {x^1}\geq 0\nonumber\\
& &\qquad T^1(\xi^{1}){x^1} + W^2(\xi^{1}) {x^2}(\xi^{1})= h^2(\xi^{1}),\quad \quad {x^2}(\xi^{1})\geq 0,\   \quad	\forall \xi^1\in\Xi^1 \ ,\nonumber 
\end{eqnarray}
\noindent where  $c^1\in\mathbb{R}^{n_1}$ and $h^1\in \mathbb{R}^{m_1}$ are known vectors and  $A\in \mathbb{R}^{m_1\times n_1}$ is a given (known) matrix.  
 The uncertain parameters vectors and matrices affected by the random process $\xi^1$ are then  given by
$h^{2}\in \mathbb{R}^{m_2}$, $c^2\in\mathbb{R}^{n_2}$, $T^{1}\in\mathbb{R}^{m_{2}\times n_{1}}$,
and $W^2\in \mathbb{R}^{m_2\times n_2}$.  

The goal is to find a sequence of decisions $({x^1},{x^2}(\xi^{1}))$ that minimizes the cost function in the worst-case realization of $\xi^{1}\in\Xi^1$.
The decision ${x^1}$ is selected at time 1, before the future outcome of $\xi^{1}$ has been revealed. The decision  $x^{2}$ at stage  $t=2$ is  $\mathscr{A}^{1}$-measurable and it depends only on the information up to time $2$.

We first remark that problem (\ref{RO2}) can  equivalently be rewritten as follows
\[
\label{RO2_a}
\textrm{RO}_2 =\min_{x^1} \left\{ c^{{1}^{\top}}{x^1}  + \mathscr{Q}({x^1})\text{ s.t. }  A{x^1}=h^1, \, {x^1}\geq 0\right\},
\]   
where $\mathscr{Q}$ is the \textit{worst-case recourse function}  
\[
\mathscr{Q}({x^1}):=\sup_{\xi^{1}\in\Xi^1} \mathcal{Q}\left({x^1},\xi^{1}\right),
\label{worstcaserecourse}
\]
being $\mathcal{Q}\left({x^1},\xi^{1}\right)$ the (uncertain) \textit{recourse function}
\begin{equation}
\mathcal{Q}\left({x^1},\xi^{1}\right) := \min_{{x^2}(\xi^{1})}\left\{ c^{{2}^{\top}}\!\left(\xi^{1}\right) {x^2}(\xi^{1})\text{ s.t. }   T^1(\xi^{1}){x^1} + W^2(\xi^{1}) x^2(\xi^{1})= h^2(\xi^{1}),\  {x^2}(\xi^{1})\geq 0 \right\}.\nonumber
\end{equation}

Our key observation is that problem $\textrm{RO}_2$ can be restated in the form of a so-called \textit{robust with certificates} $\textrm{RwC}_2$ problem, where we distinguish between \textit{design variables} ($x^1, \gamma$) and certificates $x^2(\xi^1)$. This observation, which represents a first result of the paper, is crucial for our successive developments and it is proved in the following Theorem.

\begin{theorem}
\label{lemma1}
The robust two-stage linear program \textrm{RO$_2$} is equivalent to the following \textit{robust with certificates} $\textrm{RwC}_2$ problem
\begin{eqnarray*}
\label{RwC2}
 \textrm{RwC}_2 &:=&\!  \min_{{x^1},\gamma}\gamma \\
&&\textrm{s.t. }  A{x^1}=h^1,\, {x^1}\ge 0\nonumber\\
& & \qquad \forall \xi^{1}\in\Xi^1,\, \exists  {x^2}(\xi^{1})\in\mathbb{R}^{n_2}\quad \textrm{\rm satisfying}\nonumber\\
& & \qquad c^{{1}^{\top}}{x^1} + c^{{2}^{\top}}\left(\xi^{1}\right) {x^2}(\xi^{1})\le \gamma \nonumber\\
& &	\qquad {x^2}(\xi^{1})\geq 0,\ T^1(\xi^{1}){x^1} + W^2(\xi^{1}) {x^2}(\xi^{1})= h^2(\xi^{1}). \nonumber
\end{eqnarray*}
\end{theorem}

\proof
We first note that Problem {RO$_2$} can be rewritten in epigraph form, by introducing the additional variable $\gamma$, as follows
\begin{eqnarray*}
\label{RO2_b}
\textrm{RO}_2& = & \! \min_{{x^1} ,\gamma}  \gamma \\
& \textrm{s.t. }&A{x^1} =h^1,\quad {x^1} \geq 0\nonumber\\
& &c^{{1}^{\top}}{x^1}+\left[ 
\begin{array}{ll}
\displaystyle \min_{ {x^2}(\xi^{1})}&c^{{2}^{\top}}\left(\xi^{1}\right) {x^2}(\xi^{1})\\
\text{s.t.}&  {x^2}(\xi^{1})\geq 0,\ T^1(\xi^{1}){x^1} + W^2(\xi^{1}) {x^2}(\xi^{1})= h^2(\xi^{1}) 
\end{array} \right]\le \gamma,\quad \forall \xi^1\in\Xi^1, 
\nonumber 
\end{eqnarray*}
or, noting that $c^{{1}^{\top}}{x^1}$ does not depend on $\xi^1$,  as
\begin{eqnarray*}
\label{RO2_c}
\textrm{RO}_2& = & \! \min_{{x^1} ,\gamma}  \gamma \\
& \textrm{s.t. }&A{x^1} =h^1,\quad {x^1} \geq 0\nonumber\\
& &\left[ 
\begin{array}{ll}
\displaystyle \min_{ {x^2}(\xi^{1})}&c^{{1}^{\top}}{x^1}+c^{{2}^{\top}}\left(\xi^{1}\right) {x^2}(\xi^{1})\\
\text{s.t.}&  {x^2}(\xi^{1})\geq 0,\ T^1(\xi^{1}){x^1} + W^2(\xi^{1}) {x^2}(\xi^{1})= h^2(\xi^{1}) 
\end{array} \right]\le \gamma,\quad \forall \xi^1\in\Xi^1,
\nonumber 
\end{eqnarray*}
or, equivalently, as
\begin{eqnarray*}
\label{RO2_X}
 \textrm{RO}_2 &=&\!  \min_{{x^1},\gamma}\gamma  \\
& &\textrm{s.t. }  A {x^1}=h^1,\quad {x^1} \geq 0 \nonumber\\
& &	\qquad ({x^1},\gamma) \in \mathcal{X}_{\mathrm{RO}_2}(\xi^{1}),\  \forall \xi^{1} \in \Xi^1, \nonumber
\end{eqnarray*}
where the set $\mathcal{X}_{\mathrm{RO}_2}(\xi^{1})$ is defined as
\begin{displaymath}
\mathcal{X}_{\mathrm{RO}_2}(\xi^{1}):=\left\{ 
({x^1},\gamma)\in\mathbb{R}^{n_1+1}_+ \text{ s.t. }
\left[
\begin{array}{ll}
\displaystyle \min_{ {x^2}(\xi^{1})}&c^{{1}^{\top}}{x^1}+c^{{2}^{\top}}\left(\xi^{1}\right) {x^2}(\xi^{1})\\
\text{s.t.}&  {x^2}(\xi^{1})\geq 0,\
T^1(\xi^{1}){x^1} + W^2(\xi^{1}) {x^2}(\xi^{1})= h^2(\xi^{1}) 
\end{array} \right]\le \gamma
\right\}\ .
\end{displaymath}
Similarly, Problem RwC$_2$ rewrites
\begin{eqnarray*}
\label{RwC_X}
 \textrm{RwC}_2 &=&\!  \min_{{x^1},\gamma}\gamma  \\
& &\textrm{s.t. }  A {x^1}=h^1,\quad {x^1} \geq 0 \nonumber\\
& &	\qquad ({x^1},\gamma) \in \mathcal{X}_{\mathrm{RwC}_2}(\xi^{1}),\  \forall \xi^{1} \in \Xi^1, \nonumber
\end{eqnarray*}
where the set $\mathcal{X}_{\mathrm{RwC}_2}(\xi^{1})$ is defined as
\begin{displaymath}
\mathcal{X}_{\mathrm{RwC}_2}(\xi^{1}):=\left\{ 
({x^1},\gamma)\in\mathbb{R}^{n_1+1}_+ \text{ s.t. }
\left\{\begin{array}{l}
 \exists {x^2}(\xi^{1})\in\mathbb{R}^{n_2}_+\textrm{ satisfying  }\nonumber \\ 
 c^{{1}^{\top}}{x^1} + c^{{2}^{\top}}\left(\xi^{1}\right) {x^2}(\xi^{1})\leq \gamma\nonumber\\
 T^1(\xi^{1}){x^1} + W^2(\xi^{1}){x^2}(\xi^{1})= h^2(\xi^{1}) 
\end{array} 
\right.
\right\} \ .
\end{displaymath}
So, we just need to prove that $\mathcal{X}_{\textrm{RO}_2}(\xi^{1})\equiv \mathcal{X}_{\textrm{RwC}_2}(\xi^{1})$ for the minimum value of $\gamma$.
\noindent\begin{itemize}
\item[$\bullet$] We prove that if $({x^1},\gamma)\in \mathcal{X}_{\mathrm{RO_2}}$, then $({x^1},\gamma)\in \mathcal{X}_{\mathrm{RwC_2}}$. 
If $({x^1},\gamma)\in \mathcal{X}_{\textrm{RO}_2}$, then $\exists x^{2}(\xi^1)\in \mathbb{R}_{+}^{n_2}$ such that $T^1(\xi^{1}){x^1} + W^2(\xi^{1}){x^2}(\xi^{1})= h^2(\xi^{1})$ is satisfied and   
$$\min_{ {x^2}(\xi^{1})}c^{{1}^{\top}}{x^1}+c^{{2}^{\top}}\left(\xi^{1}\right) {x^2}(\xi^{1})\leq \gamma,$$ for the minimum value of $\gamma$. 
Consequently $({x^1},\gamma)\in \mathcal{X}_{\textrm{RwC}_2}$.
\item[$\bullet$] Conversely if $({x^1},\gamma)\in \mathcal{X}_{\mathrm{RwC}_2}$, then we need to prove that $({x^1},\gamma)\in \mathcal{X}_{\mathrm{RO}_2}$.
If $({x^1},\gamma)\in \mathcal{X}_{\textrm{RwC}_2}$ then $\exists x^{2}(\xi^1)\in \mathbb{R}_{+}^{n_2}$ such that $T^1(\xi^{1}){x^1} + W^2(\xi^{1}){x^2}(\xi^{1})= h^2(\xi^{1})$ and $c^{{1}^{\top}}{x^1}+c^{{2}^{\top}}\left(\xi^{1}\right) {x^2}(\xi^{1})\leq \gamma $ for the minimum value of $\gamma$. This implies that ${x^2}(\xi^{1})$ is the minimum of   $c^{{1}^{\top}}{x^1}+c^{{2}^{\top}}\left(\xi^{1}\right) {x^2}(\xi^{1})$. By contradiction if ${x^2}(\xi^{1})$ were not be the minimum then $\gamma$ would not be at the minimum of problem $\textrm{RwC}_2$. 
\end{itemize}
\qed

Based on the result of Theorem  \ref{lemma1}, we are now ready to formulate the scenario with certificates counterpart of problem $\textrm{RO}_2$. 
To this end, we exploit the probabilistic information about the uncertainty and, similarly to what proposed in  \cite{Vayanos2012459}, we 
adopt a sampling approach, based on the extraction of $N$ independent identically distributed (iid) samples 
\[
\xi^{{1}^{(1)}},\dots,\xi^{{1}^{(N)}}
\]
 of the  random variable   $\xi^{1}$.
The  samples are extracted according to the probability distribution of the uncertainty over $\Xi^1$. 
Let $T^1(\xi^{{1}^{(i)}})$, $h^2(\xi^{{1}^{(i)}})$, $c^2(\xi^{{1}^{(i)}})$ be the realization of $T^1(\xi^1)$, $h^2(\xi^1)$ and $c^2(\xi^1)$ under scenario $\xi^{{1}^{(i)}}$, $i=1,\dots,N$,  and let $x^2_i$ be the certificate variables  created for the samples $\xi^{{1}^{(i)}}$, $i=1,\dots,N$.
These samples are used to construct the following \textit{scenario with certificates} $\textrm{SwC}^N_2$ problem based on $N$ instances (scenarios) of the uncertain constraints
\begin{eqnarray}
\label{sp_twoNbis}
 \textrm{SwC}^N_2 &:=&\! \min_{{x^1},\gamma,x^2_1,\dots, x^2_N} \gamma \\
& & \textrm{s.t. }			 A{x^1}=h^1,\quad {x^1}\geq 0\nonumber\\
& & \qquad c^{{1}^{\top}}{x^1} + c^{{2}^{\top}}\left(\xi^{{1}^{(i)}}\right) x^2_i\leq \gamma	\nonumber\\ 
& &\qquad T^1(\xi^{{1}^{(i)}}){x^1} + W^2(\xi^{{1}^{(i)}}) x^2_i= h^2(\xi^{{1}^{(i)}}),\quad \quad x^2_i\geq 0,\quad	i=1,\dots,N. \nonumber
\end{eqnarray}
The solution of problem $\textrm{SwC}^N_2$ is denoted with $(\hat{x}_N^1,\hat{\gamma}_N)$.
We note that in problem $\textrm{SwC}^N_2$, a different \textit{certificate} $x_2^i$ is constructed for any sample $\xi^{{1}^{(i)}}$. The rationale behind this approach is the following: we are not interested in the explicit knowledge of the function $x^{2}(\xi^1)$, what we are content with is that for every possible value of the uncertainty \textit{there exists} a possible choice of $x_2$ compatible with the ensuing realization of the constraints. In the SwC approach, this requirement is relaxed by asking that this is true only for the sampled scenario.
Note that this represents a key difference with respect to other sampling based approaches.
In particular, in \cite{Vayanos2012459} different explicit parameterizations of the function $x^{2}(\xi^1)$ are introduced, of the form 
\[
x^{2}(\xi^1)= \sum_{k=1}^{M} c_k \phi^2_k(\xi^1),
\]
where $\phi^2_1,\dots,\phi_M^2$ are specific basis functions, which can be for instance algebraic polynomials, trigonometric polynomials, sigmoidal  or gaussian radial basis functions and  $c_k$ represent the coefficients of the linear combinations, which become the new decision variables. 
It is easy to infer how this latter approach is bound of being more conservative, since the existence of a solution with a pre-specified form is required.

It is clear that the  approximate solution returned by problem $\textrm{SwC}_2^N$ is optimistic, since it considers only a subset of possible scenarios. That is, the following bound holds for all $N$:
\begin{equation}
\label{SwC_lowerbound}
\textrm{SwC}_2^N \le \textrm{RO}_2.
\end{equation}
Hence, we have derived a lower bound, which by construction is better than bounds derived using wait-and-see approaches, as discussed in Section \ref{boundsMultistageRO}.
{Moreover, it is easy to show that the formulation is consistent, that is 
\[
\lim_{N\to\infty} \textrm{SwC}_2^N = \textrm{RO}_2.
\]
}
More importantly, we note that, by exploiting recent results in \cite{Formentin2015RobustSA}, it is possible to provide a formal assessment about its probabilistic properties. To this end, let formally introduce  the  violation probability  $\textrm{V}_2({x^1},\gamma)$ of $(x^1,\gamma)$ as follows

\begin{displaymath}
\textrm{V}_2({x^1},\gamma) :=\pr\left\{ 
\exists\xi^1\in\Xi^1 \text{ for which } \nexists {x^2}(\xi^{1})\in\mathbb{R}^{n_2}_+ :
\left[\begin{array}{l}
 c^{{1}^{\top}}{x^1} + c^{{2}^{\top}}\left(\xi^{1}\right) {x^2}(\xi^{1})\leq \gamma\nonumber\\
 T^1(\xi^{1}){x^1} + W^2(\xi^{1}){x^2}(\xi^{1})= h^2(\xi^{1}) 
\end{array} 
\right.
\right\}.
\end{displaymath}

The interpretation of the violation probability of the solution $x^1$ is as follows: if we select as first stage solution $x^1$, then $\textrm{V}_2(x^1,\gamma)$ is the probability that at stage two we encounter an uncertainty realization $\xi^1$ for which there does not exist a feasible recourse decision $x^2(\xi^1)$. Clearly, the smaller is $\textrm{V}_2(x^1,\gamma)$, the higher is the probability that the solution at stage one will lead to a feasible stage two problem.  
We are in the position of providing a rigorous result connecting the violation probability to the number of samples $N$ adopted in the construction of the $\textrm{SwC}_2^N$ problem.
The following theorem holds.
 \begin{theorem}[two-stage robust linear case]
\label{the2stagesamplecomplexity}
Assume that, for any multisample extraction, the problem $\textrm{SwC}^N_2$ is feasible and attains a unique optimal solution. 
Then, given an accuracy level $\epsilon\in(0,1)$, the solution $(\hat{x}_N^1,\hat{\gamma}_N)$ of the problem (\ref{sp_twoNbis}) satisfies
\begin{equation}
\label{bound2}
\pr\left\{\textrm{V}_2({\hat{x}^1},\hat{\gamma})_{\textrm{SwC}^N_2}>\epsilon\right\}\leq B(N,\epsilon,n_1+1),
\end{equation}
where $B(N,\epsilon,n_1+1):=\sum_{k=0}^{n_1}\binom{N}{k}\epsilon^k(1-\epsilon)^{N-k}$.
\end{theorem}
\vskip 4mm

The proof of Theorem \ref{the2stagesamplecomplexity} 
follows the same lines of the results presented in 
\cite{Formentin2015RobustSA}, and 
is reported in the Appendix. The theorem provides a way to a priori bound the probability of violation of the solution of $\textrm{SwC}^N_2$.
We remark that, in the literature, the minimum number of samples for which
\eqref{bound2} holds for given 
 $\epsilon\in(0,1)$ and $\beta\in(0,1)$ is referred to as sample complexity, see \cite{IOPORT.06102166}.
Several are the results derived in literature to bound sample complexity. In  particular, in Lemma 1 and 2 in \cite{Alamo2015160}, it is proved that given $\epsilon\in(0,1)$ and $\beta\in(0,1)$
\begin{equation}
\label{Al}
N(\epsilon,\beta)\geq \frac{1}{\epsilon}{\frac{e}{e-1}}\left(\ln \frac{1}{\beta} + n_0+1\right),
\end{equation}
 where $e$ is the Euler constant.
This bound represents a (numerically) significant improvement upon other  bounds available in the literature \cite{doi:10.1137/090773490,Calafiore20111279}.

It is important to highlight that the number of samples $N$ in formula (\ref{Al}) depends only on the dimension of non-adjustable variables (or design variables); thus it reduces the number of samples needed to satisfy a prescribed level of violation with respect to that proposed in \cite{Vayanos2012459}. Indeed, in the proof of Corollary 1 in \cite{Vayanos2012459},  $N$ depends on the size of the basis $M$ and on the number of decision variables at each stage.

The results presented in this section can be readily extended to 
the more general case of dynamic  multistage ($H$-stages) robust linear decision problem under uncertainty. This is done in the next section.

\subsection{Multistage robust linear case}
\label{sec:MultistageRobustLinearCase}
We consider the following robust linear program over $H$ stages
\begin{eqnarray}
\label{ro_multi}
\textrm{RO}_H & := &  \min_{x^1,\dots,x^H(\xxi^{H-1})} \sup_{\xxi^{H-1}}  z\left[\left(x^1,\dots,x^H(\xxi^{H-1})\right),\xxi^{H-1}\right] \\
& & = \! \min_{x^1} c^1{^{\top}} \! x^1 + \nonumber \\
& & \! + \! \sup_{\xi^1\in\Xi^1}\!\!\left[\min_{x^2(\xi^1)} c^2{^{\top}}\!\!\left({\xi}^1\right) x^2\!\left({\xi}^1\right) + \sup_{\xi^2\in\Xi^2}\!\left[\dots + \!\!\!\!
										\sup_{ \xi^{H-1}\in\Xi^{H-1}}\!\!\left[\min_{x^H\!\left(\xxi^{H-1}\!\right)} {c^H}^{\top}\!\!\left({\xxi}^{H-1}\right) x^H\!\left(\xxi^{H-1}\!\right)\!\right]\!\right]\!\right]  \nonumber \\
& & \textrm{s.t. }			 A x^1=h^1,\ x^{1}\geq 0\ \nonumber\\
& &\qquad T^1({\xi}^1)x^1 + W^2(\xi^1) x^2({\xi}^1)= h^2({\xi}^1),\ \forall\xi^1\in\Xi^1\nonumber\\
& & \qquad \qquad \qquad \vdots \nonumber\\
& & \qquad  T^{H-1}(\xxi^{H-1}) x^{H-1}(\xxi^{H-2}) + W^H(\xxi^{H-1}) x^H(\xxi^{H-1})= h^H(\xxi^{H-1}),\ \forall \xxi^{H-1}\in\Xi\nonumber\\
& & \qquad  x^{t}(\xxi^{t-1})\geq 0\ ,\  t=2,\dots,H\nonumber, \quad	\forall \xxi^{t-1}\in {\sf X}_{\tau=1}^{t-1} \Xi^{\tau},
\end{eqnarray} 
\noindent where  $c^1\in\mathbb{R}^{n_1}$ and $h^1\in \mathbb{R}^{m_1}$ are known vectors and  $A\in \mathbb{R}^{m_1\times n_1}$ is  known matrix.  
The uncertain parameter vectors and matrices affected by the random process $\xi^t$ are then  given by
$h^{t}\in \mathbb{R}^{m_t}$, $c^t\in\mathbb{R}^{n_t}$, $T^{t-1}\in\mathbb{R}^{m_{t}\times n_{t-1}}$, and $W^t\in\mathbb{R}^{m_t\times n_t}$, $t=2,\dots,H$. \\

The aim of the problem $\textrm{RO}_H$  is to find a sequence of decisions $(x^1,\dots,x^H)$ that minimizes a cost function in the worst-case realization of $\xxi^{H-1}\in\Xi={\sf X}_{t=1}^{H-1} \Xi^t$. The decision process is nonanticipative and depends on the information up to time $t$ as described in Section \ref{notation}.

Similarly to the two-stage case, we first rewrite problem (\ref{ro_multi}) as the  multistage robust optimization problem with certificates $\textrm{RwC}_H$, where we distinguish between \textit{design variables} $x^1,\gamma$ and {nonanticipative} certificates $(x^2(\xi^1),\dots,x^t(\xxi^{t-1}),\dots,x^H(\xxi^{H-1}))$ as follows 
\begin{eqnarray}
\label{rwc_multi}
\textrm{RwC}_H & := &\! \min_{x^1,\gamma}\gamma 	\\
& & \textrm{s.t. }  \forall \xxi^{H-1} \in \Xi,\ \exists x^{t}(\xxi^{t-1}) \in\mathbb{R}^{n_t}_+,\  t=2,\dots,H \textrm{ satisfying  }\nonumber\\
& &	\qquad  {c^1}^{\top} x^1+{c^{2}}^\top\!\!\!\left(\xi^{1}\right) x^{2}\!\left(\xi^{1}\right)+\cdots+
{c^{H}}^\top\!\!\!\left(\xxi^{H-1}\right) x^{H}\!\left(\xxi^{H-1}\right) \le \gamma \nonumber\\
& &	\qquad  A x^1=h^1,\ x^{1}\geq 0 \nonumber\\
& &\qquad T^1({\xi}^1)x^1 + W^2(\xi^{1}) x^2({\xi}^1)= h^2({\xi}^1)\nonumber\\
& & \qquad \qquad \qquad \vdots \nonumber\\
& & \qquad  T^{H-1}(\xxi^{H-1}) x^{H-1}(\xxi^{H-2}) + W^H(\xxi^{H-1}) x^H(\xxi^{H-1})= h^H(\xxi^{H-1}).\nonumber
\end{eqnarray}

The equivalence of problems  $\textrm{RO}_H$ and $\textrm{RwC}_H$ is formally stated in the following theorem, which represents a generalization of Theorem \ref{lemma1} to the multistage case. The proof follows the same lines and is reported in Appendix B.

\begin{theorem}
\label{theomultistageequiv}
The robust multistage linear program \textrm{RO$_H$} (\ref{ro_multi}) is equivalent to the robust with certificates $\textrm{RwC}_H$ (\ref{rwc_multi}) problem.
\end{theorem}

\vskip 3mm


Again, the previous theorem is very important in that it allows to reformulate problem \textrm{RO$_H$} using  the scenario with certificates approach. 
For this purpose, we extract $N$ iid samples ${\xxi^{H-1}}^{(1)},\dots,{{\xxi^{H-1}}}^{(N)}$  according to the probability distribution of the uncertainty over $\Xi$, where 
\[
{\xxi^{H-1}}^{(i)}=({\xi^1}^{(i)},\dots,{\xi^{H-1}}^{(i)}),\ i=1,\dots,N.
\]
Let $T^{t-1}(\xxi^{{t-1}^{(i)}})$, $h^t(\xxi^{{t-1}^{(i)}})$, $c^t(\xi^{{t-1}^{(i)}})$ be the realization of $T^{t-1}(\xxi^{t-1})$, $h^{t}(\xxi^{t-1})$ and $c^t(\xxi^{t-1})$ under scenario $\xxi^{{t-1}^{(i)}}$, $i=1,\dots,N$, $t=2,\dots,H$, and 
let $x_i^t$ be the certificate  $x^t({\xxi^{t-1}}^{(i)})$ created for the sample ${\xxi^{H-1}}^{(i)}$, $i=1,\dots,N$ taking into account the history of the process until period $t-1$. That is 
\[
x_i^t=x^t({\xxi^{t-1}}^{(i)}), \quad t=2,\dots,H,
\]
which means that the decision process is still nonanticipative. 
These samples are used to construct the following \textit{multistage scenario with certificates} $\textrm{SwC}^N_H$ problem based on $N$ instances (scenarios) of the uncertain constraints
\begin{eqnarray}
\label{MSwC_N}
\textrm{SwC}_H^N& := &\! \min_{x^1,\gamma,x_i^2,\dots,x_i^H} \gamma    \\
& & \textrm{s.t. }   {c^1}^{\top} x^1+{c^2}^{\top}\left({\xi^{1}}^{(i)}\right) x^2_i  + \cdots+{c^H}^{\top}\left({\xxi^{H-1}}^{(i)}\right) x^H_i \leq \gamma,\  i=1,\dots,N \nonumber \\
& & \qquad  A x^1=h^1,\ x^{1}\geq 0 \nonumber\\
& &\qquad T^1({{\xi}^1}^{(i)})x^1 + W^2({\xi^1}^{(i)}) x^2_i= h^2({{\xi}^1}^{(i)}),\  i=1,\dots,N \nonumber\\
& & \qquad \qquad \qquad \vdots \nonumber\\
& & \qquad  T^{H-1}({\xxi^{H-1}}^{(i)}) x^{H-1}_i + W^H({\xxi^{H-1}}^{(i)}) x^H_i= h^H({\xxi^{H-1}}^{(i)}),\  i=1,\dots,N \nonumber\\
& & \qquad  x^{t}_i\geq 0\ ,\  t=2,\dots,H,\   i=1,\dots,N\nonumber.
\end{eqnarray}
The solution of problem $\textrm{SwC}_H^N$ is denoted with $(\hat{x}_N^1,\hat{\gamma}_N)$. 

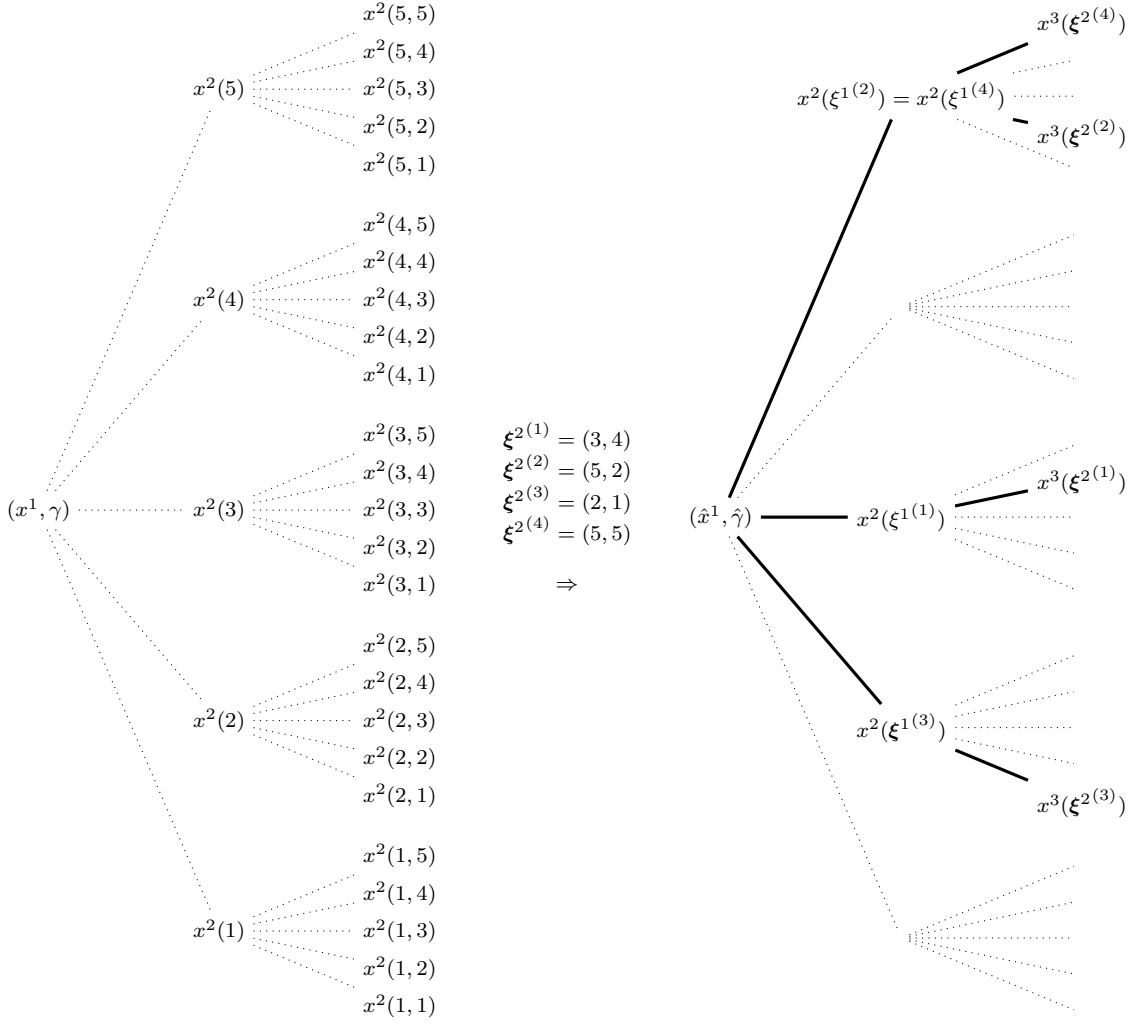
\begin{figure}
\centering
\begin{minipage}{.4\textwidth}
\begin{tikzpicture}[rotate=90, level distance=2.4cm,
  level 1/.style={sibling distance=2.8cm},
  level 2/.style={sibling distance=0.5cm}]
\tikzstyle{every node}=[]
  \node{$(x^1,\gamma)$}
    child[dotted]{node{$x^2(1)$}
      	child[dotted]{node{$x^2(1,1)$}}
      	child[dotted]{node{$x^2(1,2)$}}
      	child[dotted]{node{$x^2(1,3)$}}
      	child[dotted]{node{$x^2(1,4)$}}
      	child[dotted]{node{$x^2(1,5)$}}
    }
    child[dotted]{node{$x^2(2)$}
      	child[dotted]{node{$x^2(2,1)$}}
      	child[dotted]{node{$x^2(2,2)$}}
      	child[dotted]{node{$x^2(2,3)$}}
      	child[dotted]{node{$x^2(2,4)$}}
      	child[dotted]{node{$x^2(2,5)$}}
    }
    child[dotted]{node{$x^2(3)$}
      	child[dotted]{node{$x^2(3,1)$}}
      	child[dotted]{node{$x^2(3,2)$}}
      	child[dotted]{node{$x^2(3,3)$}}
      	child[dotted]{node{$x^2(3,4)$}}
      	child[dotted]{node{$x^2(3,5)$}}
    }
    child[dotted]{node{$x^2(4)$}
      	child[dotted]{node{$x^2(4,1)$}}
      	child[dotted]{node{$x^2(4,2)$}}
      	child[dotted]{node{$x^2(4,3)$}}
      	child[dotted]{node{$x^2(4,4)$}}
      	child[dotted]{node{$x^2(4,5)$}}
    }
    child[dotted]{node{$x^2(5)$}
      	child[dotted]{node{$x^2(5,1)$}}
      	child[dotted]{node{$x^2(5,2)$}}
      	child[dotted]{node{$x^2(5,3)$}}
      	child[dotted]{node{$x^2(5,4)$}}
      	child[dotted]{node{$x^2(5,5)$}}
    };
\end{tikzpicture}
\end{minipage}%
\begin{minipage}{.2\textwidth}$\displaystyle
\qquad
\begin{array}{c}
{\xxi^{2}}^{(1)}=(3,4)\\
{\xxi^{2}}^{(2)}=(5,2)\\
{\xxi^{2}}^{(3)}=(2,1)\\
{\xxi^{2}}^{(4)}=(5,5)\\
~\\
\Huge{\Rightarrow}
\end{array}
$
\end{minipage}%
\begin{minipage}{.4\textwidth}
\begin{tikzpicture}[rotate=90, level distance=2.4cm,
  level 1/.style={sibling distance=2.8cm},
  level 2/.style={sibling distance=0.5cm}]
\tikzstyle{every node}=[]
  \node{$({\hat{x}^1,\hat{\gamma}})$}
      child[dotted]{node{} 
      	child[dotted]{node{}}
      	child[dotted]{node{}}
      	child[dotted]{node{}}
      	child[dotted]{node{}}
      	child[dotted]{node{}}
    }
      child[very thick]{node{$x^2({\xxi^{1}}^{(3)})$}
      		child[very thick]{node{$x^3({\xxi^{2}}^{(3)})$}}
      child[thin,dotted]{node{}}
      child[thin,dotted]{node{}}
			child[thin,dotted]{node{}}
      child[thin,dotted]{node{}}
    }
    child[very thick]{node{$x^2({\xi^{1}}^{(1)})$}
      child[thin,dotted]{node{}}
      child[thin,dotted]{node{}}
			child[thin,dotted]{node{}}
			child[very thick]{node{$x^3({\xxi^{2}}^{(1)})$}}
      child[thin,dotted]{node{}}
    }
		child[dotted]{node{}
      child[dotted]{node{}}
      child[dotted]{node{}}
			child[dotted]{node{}}
			child[dotted]{node{}}
      child[dotted]{node{}}
    }
		child[very thick]{node{$x^2({\xi^{1}}^{(2)})=x^2({\xi^{1}}^{(4)})$}
      child[thin,dotted]{node{}}
      child[very thick]{node{$x^3({\xxi^{2}}^{(2)})$}}
			child[thin,dotted]{node{}}
			child[thin,dotted]{node{}}
      child[very thick]{node{$x^3({\xxi^{2}}^{(4)})$}}
    };
\end{tikzpicture}
\end{minipage}%
\caption{Example of three-stage  robust optimization problem solved through a scenario with certificates approach. In this case, the first and second {period} uncertainties $\xi^1$ and $\xi^2$ can assume the values $\{1,2,3,4,5\}$, with equal probability.  On the left, the complete (robust) tree for problem $\textrm{RO}_3$ is shown. On the right, the  $\textrm{SwC}_3^4$, based on the extraction of $N=4$ samples (thick lines) of the uncertain constraints in the initial problem {is shown}. In the example, the sampled uncertainties extracted are ${\xxi^{2}}^{(1)}=(3,4)$; $
{\xxi^{2}}^{(2)}=(5,2)$; $
{\xxi^{2}}^{(3)}=(2,1)$; $
{\xxi^{2}}^{(4)}=(5,5)$. The quantities $x^2({\xi^1}^{(i)})$, $x^3({\xxi^{2}}^{(i)})$ represent the certificates over the  samples $i=1,2,3,4$. 
Notice that the extracted samples ${\xi^{1}}^{(2)}={\xi^{1}}^{(4)}=5$ coincide, and 
 the  scenario with certificates $\textrm{SwC}_3^4$ is constructed accordingly so to satisfy the   non-anticipativity constraint  $x^2({\xi^{1}}^{(2)})=x^2({\xi^{1}}^{(4)})$.}
\label{FigureTree}
\end{figure}

\begin{remark}[Scenario construction] \it Note that the type of scenario construction proposed by the implementation of problem $\textrm{SwC}_H^N$ differs from the classical \textit{scenario trees} proposed in literature. Indeed, instead of generating a few possible ``leaves''at every stage, and considering the tree obtained from all possible combinations, we sample $N$ different ``paths''. 
This procedure is illustrated in Figure \ref{FigureTree}, which shows the construction of $\textrm{SwC}_H^N$ from $\textrm{RwC}_H$ in the case of a three-stage  robust optimization problem in which the first and second {period} uncertainties {are discrete and} can take a finite number of possible values. This allows to visualize the tree of all possible solution (left figure). The figure on the right shows the paths generated by a scenario with certificates $\textrm{SwC}_3^4$,  based on $N=4$ samples (thick lines) of the uncertain constraints in the initial problem.  $(\hat{x}^1,\hat{\gamma})$  represent the design variables solution of $\textrm{SwC}_3^4$,  and  $(x^2({\xi^1}^{(i)}),x^3({\xxi^{2}}^{(i)}))$ the certificates over the  samples $i=1,2,3,4$.  Notice that, the  nonanticipativity constraints 	{have} to be imposed, which in our case translate in requiring that  $x^2({\xi^{1}}^{(2)})=x^2({\xi^{1}}^{(4)}).$
\end{remark}

Again, by construction, the following bounds hold
\begin{equation}
\label{SwC_lowerbound_H}
\textrm{SwC}_H^{N_1} \le \textrm{SwC}_H^{N_2} \le \textrm{RO}_H, \quad 1\le N_1\le N_2,
\end{equation}
where we explicitly highlight that the lower bound improves for increasing values of $N$.
{In particular, it can be shown that  
\[
\lim_{N\to\infty} \textrm{SwC}_H^N = \textrm{RO}_H.
\]
}
Moreover, similarly to the two-stage case, we can formally investigate the probabilistic properties of the approximate solution returned by problem $\textrm{SwC}_H^N$. To this end, we introduce  the reliability $\textrm{R}_H({x^1},\gamma)$ and violation probability of the scenario with certificates problem as follows\\

$
\textrm{V}_H(x^1,\gamma)  :=  
$
\[
\pr\left\{
\begin{array}{ll} 
\exists\xxi^{H-1}\in \Xi \text{ for which }&
 \nexists\  x^{t}(\xxi^{t-1}) \in\mathbb{R}^{n_t}_+,\ t=2,\dots,H :\\
&
\left[
\begin{array}{l}
 {c^1}^{\top} x^1+{c^{2}}^\top\!\!\!\left(\xi^{1}\right) x^{2}\!\left(\xi^{1}\right)+\cdots+
{c^{H}}^\top\!\!\!\left(\xxi^{H-1}\right) x^{H}\!\left(\xxi^{H-1}\right) \le \gamma \nonumber\\
 T^{t-1}(\xxi^{t-1})x^{t-1}(\xxi^{t-2}) + W^t(\xxi^{t-1}) x^{t}(\xxi^{t-1})= h^t(\xxi^{t-1}),\\
t=2,\dots,H  
\end{array} 
\right.
\end{array}
\right\}.  
\]

We provide now a sample complexity result for the multistage robust linear case which extends Theorem \ref{the2stagesamplecomplexity} for the two-stage robust linear case. The proof is given in Appendix B.

\begin{theorem}[multistage robust linear case]
Assume that, for any multisample extraction, problem $\textrm{SwC}_H^N$ is feasible and attains a unique optimal solution. Then, given an accuracy level $\epsilon\in(0,1)$, the solution $(\hat{x}^1,\hat{\gamma})$ of problem (\ref{MSwC_N}) satisfies
\begin{equation*}
\pr\left\{\textrm{V}({\hat{x}^1,\hat{\gamma}})>\epsilon\right\}\leq B(N,\epsilon,n_0+1),
\end{equation*}
where $B(N,\epsilon,n_0+1):=\sum_{k=0}^{n_0}\binom{N}{k}\epsilon^k(1-\epsilon)^{N-k}$.
\label{theomulti}
\end{theorem}

We note that the sample complexity for guaranteeing with high probability $(1-\beta)$ that the solution of problem $\textrm{SwC}_H^N$ has a violation probability bounded by $\epsilon$
can be computed by \eqref{Al}. It is important to remark again that, also in the multistage case, the necessary number of samples $N$ \textit{does not depend on the number of stages}
$H$. This is in sharp contrast with the setup in \cite{Vayanos2012459}, in which $N$ depends on $\sum_{i=0}^H n_i\times M_i$, that is on the number of decision variables at each stage multiplied by the number of basis functions chosen for each stage. 
On the other hand,  problem $\textrm{SwC}_H^N$ introduces an increment in the number of variables, since new variables are introduced for each stage. This growth can be easily handled in the case of linear programs, which constitute the main focus of this paper. 
We observe however that the SwC setup can be easily extended to the general context of convex multistage problems. This is briefly outlined in the next section.

\subsection{Extension to the multistage robust convex case}
\label{sec:ConvexCase}
In this section we further generalize the formulation given in Section \ref{sec:MultistageRobustLinearCase} to a  dynamic  multistage ($H$-stages) robust convex decision  problem under uncertainty, which can be formulated as follows
\begin{eqnarray*}
\textrm{CRO}_H& := &\! \min_{x^1, x^2(\xi^{1}),\dots,x^{H}(\xxi^{H-1})} \sup_{\xxi^{H-1}\in\Xi} f(x^1,x^2(\xi^{1}),\dots,x^H(\xxi^{H-1}),\xxi^{H-1})  \\
& & \textrm{s.t. }			g(x^1,x^2(\xi^{1}),\dots,x^H(\xxi^{H-1}),\xxi^{H-1})\leq  0,\ \forall \xxi^{H-1}\in\Xi  \nonumber\\
& & \qquad x^{1}\geq 0\ ,\quad x^{t}(\xxi^{t-1})\geq 0\ ,\  t=2,\dots,H,\nonumber
\label{ro_convex}
\end{eqnarray*} 
\noindent where  $f:\mathbb{R}^{\sum_{t=1}^{H} n_t} \times \Xi \to  \mathbb{R}$ and $g: \mathbb{R}^{\sum_{t=1}^{H} n_t} \times \Xi \to  \mathbb{R}$ are convex in $x^t\in\mathbb{R}^{{n}_t}_+ $, $t=1,\dots,H$ and continuous in $(x^t,\xxi^{H-1})$. Again we assume that the decision process in nonanticipative according to the desciption given in Section \ref{notation}.

The aim of  problem $\textrm{CRO}_H$ is to find a sequence of decisions $(x^1,\dots,x^H)$ that minimizes a cost function $f$ in the worst-case realization of $\xxi^{H-1}\in\Xi={\sf X}_{t=1}^{H-1} \Xi^t$.
First, we observe that  problem $\textrm{CRO}_H$ can be rewritten as the following  convex multistage robust optimization problem with certificates ($\textrm{CRwC}_H$), where again we distinguish between \textit{design variables} $x^1,\gamma$ and certificates $(x^2(\xi^1),\dots,x^H(\xxi^{H-1}))$ as follows 
\begin{eqnarray*}
\label{crwc_multi}
\textrm{CRwC}_H & := &\! \min_{x^1,\gamma}\gamma 	\\
& & \textrm{s.t. }  \forall \xxi^{H-1} \in \Xi,\ \exists x^{t}(\xxi^{t-1}) \in\mathbb{R}^{n_t}_+,\  t=2,\dots,H \textrm{ satisfying  }\nonumber\\
& &	\qquad  f(x^1,x^2(\xi^{1}),\dots,x^H(\xxi^{H-1}),\xxi^{H-1}) \le \gamma \nonumber\\
& &\qquad 	g(x^1,x^2(\xi^{1}),\dots,x^H(\xxi^{H-1}),\xxi^{H-1})\leq  0\nonumber\\
& & \qquad x^{1}\geq 0\ ,\quad x^{t}(\xxi^{t-1})\geq 0\ ,\  t=2,\dots,H\nonumber.
\end{eqnarray*}
Then, we extract $N$ iid samples ${\xxi^{H-1}}^{(1)},\dots,{{\xxi^{H-1}}}^{(N)}$, and denote by $x_i^t$  the certificate  $x^t({\xxi^{t-1}}^{(i)})$ created for  sample ${\xxi^{H-1}}^{(i)}$, $i=1,\dots,N$. This means that the decision process is still nonanticipative.
These samples are used to construct the following \textit{multistage convex scenario with certificates} $\textrm{CSwC}^N_H$ problem
\begin{eqnarray}
\label{CSwC_N}
\textrm{CSwC}_H^N& := &\! \min_{x^1,\gamma,x_i^2,\dots,x_i^H} \gamma    \\
&   &  \textrm{s.t. }   f(x^1,x^2_i,\dots,x^H_i,\xxi^{(i)}) \leq \gamma,\  i=1,\dots,N \nonumber \\
& &\qquad  g(x^1,x^2_i,\dots,x^H_i,\xxi^{(i)}) \leq 0 ,\  i=1,\dots,N \nonumber\\
& & \qquad x^{1}\geq 0\ ,\quad x^{t}_i\geq 0\ ,\  t=2,\dots,H,\   i=1,\dots,N\nonumber.
\end{eqnarray}
The solution of problem $\textrm{CSwC}_H^N$ is denoted with $(\hat{x}_N^1,\hat{\gamma}_N)$.
In order to  investigate  probabilistic properties of the approximate solution returned by problem $\textrm{CRwC}_H^N$ we introduce the violation probability $\textrm{V}_H(x^1,\gamma)$ of its solution $(x^1,\gamma)$
\[
\textrm{V}_H(x^1,\gamma) :=
\pr\left\{
\begin{array}{ll} 
\exists\xxi^{H-1}\in \Xi \text{ for which }&
 \nexists\  x^{t}(\xxi^{t-1}) \in\mathbb{R}^{n_t}_+,\ t=2,\dots,H :\\
&\left[\begin{array}{l}
f(x^1,x^2(\xi^{1}),\dots,x^H(\xxi^{H-1}),\xxi^{H-1}) \le \gamma \\
g(x^1,x^2(\xi^{1}),\dots,x^H(\xxi^{H-1}),\xxi^{H-1})\leq  0  \\
t=2,\dots,H
\end{array} \right.
\end{array} \right\}.
\]
Then, the following sample complexity result for the multistage robust convex  programs holds true:

\begin{corollary}[multistage robust convex case]
Assume that, for any multisample extraction,  problem $\textrm{CSwC}_H^N$ is feasible and attains a unique optimal solution. Then, given an accuracy level $\epsilon\in(0,1)$, the solution $(\hat{x}^1,\hat{\gamma})$ of problem (\ref{CSwC_N}) satisfies
\begin{equation*}
\pr\left\{\textrm{V}_H({\hat{x}^1,\hat{\gamma}})>\epsilon\right\}\leq B(N,\epsilon,n_0+1)
\end{equation*}
where $B(N,\epsilon,n_0+1):=\sum_{k=0}^{n_0}\binom{N}{k}\epsilon^k(1-\epsilon)^{N-k}$.
\end{corollary}
\proof The proof works similarly to the one of Theorem \ref{theomulti} for the multistage robust linear case and is omitted for brevity.
 \qed

\section{Lower Bounds for Multistage Linear  Robust Optimization Problems}
\label{boundsMultistageRO}
In this section, we present the robust counterpart of different lower bounds known in the context of stochastic programming, see for instance 
 \cite{MAB2013}, \cite{CMS2016} and \cite{Maggioni2016}. To the best of our knowledge such relaxations, while  frequently encountered when facing {stochastic} multistage problems,  have never been formally stated in the context of robust programming. In particular, we here  introduce  and  compare them in terms of optimal objective function values for the case of robust multistage linear programs. Similarly lower bounds for multistage convex robust programs can be defined.
 
First, we introduce the \textit{robust multistage wait-and-see} problem $\textrm{RWS}_H$, where  the realizations  of all the random parameters are assumed to be known at the first stage, which takes the form
%
\begin{eqnarray}
\label{WS_multi_conditional_exp}
\textrm{RWS}_H 
& :=   \sup_{\xxi^{H-1}}\quad & \min_{(x^1(\xxi^{H-1}),\dots,x^H(\xxi^{H-1}))} z\left[\left(x^1(\xxi^{H-1}),\dots,x^H(\xxi^{H-1})\right),\xxi^{H-1}\right]\\
 &:=  \sup_{\xxi^{H-1}} \quad & \min  _{ x^1(\xxi^{H-1}), \dots, x^H(\xxi^{H-1})}	 {c^1}^{\top} x^1(\xxi^{H-1}) +\! \dots \!+
 {c^H}^{\top} x^H(\xxi^{H-1})\nonumber \\
 & &\textrm{s.t. }	 A x^1=h^1,\ x^{1} \geq 0\nonumber\\
& &  \qquad T^1(\xi^{1})x^1(\xxi^{H-1}) + W^2 (\xi^{1}) x^2(\xxi^{H-1})= h^2(\xi^{1}),\
 \nonumber\\
& & \qquad \vdots\nonumber\\
& & \qquad T^{H-1}(\xxi^{H-1}) x^{H-1}(\xxi^{H-1}) + W^H (\xxi^{H-1}) x^H(\xxi^{H-1})\!=\!
h^H(\xxi^{H-1})\ 
 \nonumber\\
& &\qquad  x^{t}(\xxi^{H-1})\geq 0\ ,\  t=2,\dots,H, \label{RWS}
\end{eqnarray}
{where with $z\left[\left(x^1(\xxi^{H-1}),\dots,x^H(\xxi^{H-1})\right),\xxi^{H-1}\right]$ we denote in  a compact way the objective function and constraints of problem (\ref{RWS}).}
Notice that, in the above setup, the minimum and supremum have been exchanged. Hence, the decision process is \textit{anticipative}, since the decisions $x^1,x^2,\dots,x^H$ depend on  a given realization of $\xxi^{H-1}$. We introduce the following definition, which is an immediate extension of the concept of \textit{Expected Value of Perfect Information} for stochastic programs,
\begin{definition}
The difference 
\begin{equation}
{\textrm{RVPI}_H} :=  \textrm{RO}_H - \textrm{RWS}_H,
\end{equation} 
denotes the \textit{Robust Value of Perfect Information}  and compares robust multistage wait-and-see $\textrm{RWS}_H$ and robust multistage $\textrm{RO}_H$. 
\end{definition}

The {$\textrm{RVPI}_H$} can be interpreted as a measure of the advantage of reaching perfect information: a small {$\textrm{RVPI}_H$} indicates a small advantage for reaching the perfect information since all possible realizations have similar costs. In particular, the following inequality can be proven.

\begin{proposition}[lower bound for $\textrm{RO}_H$] 
\label{teo1}
Given the robust multistage linear optimization problem $\textrm{RO}_H$ defined in (\ref{ro_multi}), and the robust multistage wait-and-see problem  $\textrm{RWS}_H$ defined in (\ref{WS_multi_conditional_exp}), the following inequality holds true
\begin{equation}
\textrm{RWS}_H\leq \textrm{RO}_H.
\label{RWSRO}
\end{equation}
\end{proposition}
\proof For every realization, $\xxi^{H-1}$, we have the relation
\begin{equation*}
z\left[\left(\tilde{x}^1(\xxi^{H-1}),\dots,\tilde{x}^H(\xxi^{H-1})\right),\xxi^{H-1}\right] \leq  z\left[\left({x^1}^{\ast},\dots,{x^H}^{\ast}(\xxi^{H-1})\right),\xxi^{H-1}\right],
\end{equation*}
where,  $\left({x^1}^{\ast},\dots,{x^H}^{\ast}(\xxi^{H-1})\right)$ denotes an optimal solution to the $RO_H$ problem (\ref{ro_multi}) and \\ $\left(\tilde{x}^1(\xxi^{H-1}),\dots,\tilde{x}^H(\xxi^{H-1})\right)$ denotes the optimal solution for each  realization of $\xxi^{H-1}$. 
Taking the supremum of both sides yields the required inequality. \qed

	\vspace{0.5cm}
A second lower bound for problem $\textrm{RO}_H$ can be obtained by relaxing the nonanticipativity constraints only in stages $2,\dots,H$ {(see \cite{MAB2013})}. The ensuing program is the so-called \textit{robust two-stage relaxation} {$\textrm{RT}_H$}. 
Formally, consider the discrete random process  as follows  
$$\tilde{\xxi}^{\ t_{-}}:=(\xi^1,\tilde{\xi}^2,\ldots,\tilde{\xi}^{t}),\quad t=2,\ldots,H-1,$$
where $\tilde{\xi}^{t}$ is a deterministic realization of the random process $\xi^{t}$.
For instance, for long processes, $\tilde{\xi}^t$, {$t=2,\ldots,H-1$} can be chosen as the expected value of the random variable $\xi^t$.
\noindent We denote the robust two-stage relaxation problem $\textrm{RT}_H$, as follows 

\begin{eqnarray*}
\label{RT}
\textrm{RT}_H \!& := &\!\min_{x^1}  {c^1}^{\top} x^1 \!+\! \sup_{\xi^1}\!\!\left[\!\min_{x^2,\ldots,x^H}\!\!\! {c^2}^{\top} x^2\!\left(\xi^1\right) + {c^3}^{\top} x^3(\tilde{\xxi}^{{2}_{-}}) +\ldots + {c^H}^{\top} x^H\!(\tilde{\xxi}^{{H-1}_{-}})\!  \right]  \\
& & \textrm{s.t. }			 A x^1=h^1,\ x^{1}\geq 0\nonumber\\
& &\qquad T^1(\xi^1)x^1 \!+\! W^2(\xi^1) x^2(\xi^1)= h^2(\xi^1),\ \forall\xi^1\in\Xi^{1} \nonumber \\
& &\qquad \vdots\nonumber\\
& &\qquad T^{H-1}\!(\tilde{\xxi}^{{H-1}_{-}}\!) x^{H-1}\!(\tilde{\xxi}^{{H-1}_{-}}\!)\! +\! W^H\! (\tilde{\xxi}^{{H-1}_{-}}\!) x^H\!(\tilde{\xxi}^{{H-1}_{-}}\!)\!\!=\!\!
h^H\!(\tilde{\xxi}^{{H-1}_{-}}\!),\forall\xi^1\in\Xi^{1}\nonumber\\
& &\qquad  x^{t}(\tilde{\xxi}^{\ {t-1}_{-}})\geq 0,\  t=2,\dots,H,\ \forall\xi^1\in\Xi^{1}.\nonumber
\end{eqnarray*}
Following reasonings similar to those in the proof {Proposition \ref{teo1}}, based on relaxation of constraints respectively in the first stage and in the following ones, the following {bounds} can be proven.

\begin{proposition}[Chain of lower bounds for $\textrm{RO}_H$]
\label{teo11}
Given the robust multistage linear optimization problem $\textrm{RO}_H$ (\ref{ro_multi}), the robust multistage wait-and-see problem  $\textrm{RWS}_H$ (\ref{WS_multi_conditional_exp}) and the robust two-stage relaxation problem $\textrm{TP}_H$,  the following inequalities hold true
\begin{equation}
\textrm{RWS}_H \leq  \textrm{RT}_H \leq \textrm{RO}_H.
\label{RWSRP}
\end{equation}
\end{proposition}

\vspace{0.5cm}

We remark that, in the general case, both problems $\textrm{RWS}_H$ and $\textrm{RT}_H$ may be hard to solve. In such case, one can recur to sampled versions of them. In particular, we can introduce the \textit{sampled wait-and-see}
problem $\textrm{SWS}_H^N$, based on the extraction of $N$ iid samples ${\xxi^{H-1}}^{(1)},\dots,{{\xxi^{H-1}}}^{(N)}$ 

\begin{eqnarray}
\label{Scenario_WS}
\textrm{SWS}_H^N
:= & \sup_{i=1,\ldots,N} \quad&\min  _{ x^1({\xxi^{H-1}}^{(i)}), \dots, x^H({\xxi^{H-1}}^{(i)})}	 {c^1}^{\top} x^1({\xxi^{H-1}}^{(i)}) +\! \dots \!+
 {c^H}^{\top} x^H({\xxi^{H-1}}^{(i)}) \\
 & &\textrm{s.t. }	 A x^1=h^1,\ x^{1} \geq 0\nonumber\\
& &  \qquad T^1({\xi^{1}}^{(i)})x^1({\xxi^{H-1}}^{(i)}) + W^2 ({\xi^{1}}^{(i)}) x^2({\xxi^{H-1}}^{(i)})= h^2({\xi^{1}}^{(i)})\nonumber\\
& & \qquad \vdots\nonumber\\
& & \qquad T^{H-1}({{\xxi^{H-1}}^{(i)}}) x^{H-1}({\xxi^{H-1}}^{(i)}) + W^H ({{\xxi^{H-1}}^{(i)}}) x^H({\xxi^{H-1}}^{(i)})\!=\!
h^H({{\xxi^{H-1}}^{(i)}})\  \nonumber\\
& &\qquad  x^{t}({\xxi^{H-1}}^{(i)})\geq 0\ ,\  t=2,\dots,H.\nonumber
\end{eqnarray}

We note that probabilistic guarantees of the solution returned by problem $\textrm{SWS}_H^N$ can be directly derived using the maximization bound in \cite{TeBaDa:97}.
Similarly, one can extract $N$ iid samples of the first {period random process} ${\xi^1}^{(1)},\dots,{\xi^1}^{(N)}$, and construct the scenario with certificates version of the $\textrm{RT}_H^N$ problem
\begin{eqnarray}
\label{SwCT}
\textrm{SwCT}_H^N \!& := &\!\min_{x^1,\gamma,x^2_i,\dots,x^H_i} \gamma  \\
& & \textrm{s.t. }  {c^1}^{\top} x^1 \!+\!  {c^2}^{\top} x^2_i + {c^3}^{\top} x^3_i +\ldots + {c^H}^{\top} x^H_i\leq \gamma,\  i=1,\dots,N \nonumber\! \\
& &\qquad 	A x^1=h^1,\ x^{1}\geq 0\nonumber\\
& &\qquad T^1({\xi^1}^{(i)})x^1 \!+\! W^2({\xi^1}^{(i)}) x^2_i= h^2({\xi^1}^{(i)}),\  i=1,\dots,N \nonumber\\
& &\qquad \vdots\nonumber\\
& &\qquad T^{H-1}\!(\tilde{\xxi}^{{H-1}_{-}(i)}\!) x^{H-1}_i +\! W^H\! (\tilde{\xxi}^{{H-1}_{-}(i)}\!) x^H_i\!\!=\!\!
h^H\!(\tilde{\xxi}^{{H-1}_{-}(i)}\!),\  i=1,\dots,N \nonumber\\
& &\qquad x^{t}_i\geq 0,\  i=1,\dots,N,\  t=2,\dots,H\nonumber,
\end{eqnarray}
{where $\tilde{\xxi}^{{H-1}_{-}(i)}$ $i=1,\dots,N$ are iid samples of $\tilde{\xxi}^{{H-1}_{-}}$.}
Again, probabilistic guarantees of the solution of problem $\textrm{SwCT}^N_H$ being also a solution of   $\textrm{RT}_H$ can be obtained on the same lines of Theorem \ref{theomulti}.  We conclude this section by providing the following proposition, which shows the relationship between the various lower bounds based on sampling presented in this paper.

\begin{proposition}[Chain of sampling-based lower bounds for $\textrm{RO}_H$]
\label{teo11bis}
Given the robust multistage linear optimization problem $\textrm{RO}_H$ (\ref{ro_multi}), the scenario with certificates relaxation $\textrm{SwC}_H^N$ (\ref{MSwC_N}), the sampled multistage wait-and-see problem $\textrm{SWS}_H^N$ \eqref{Scenario_WS},  and the scenario with certificates two-stage relaxation $\textrm{SwCT}_H^N$ \eqref{SwCT}, the following chain of inequalities holds true
\begin{equation}
\textrm{SWS}_H^N \leq  \textrm{SwCT}_H^N\le \textrm{SwC}_H^N \leq \textrm{RO}_H.
\label{RWSRPbis}
\end{equation}
\end{proposition}


\section{Numerical Results: Inventory Management with Cumulative Orders}
\label{numericalresults}
In this section, to show the effectiveness of the proposed approach, we
consider a problem from inventory management which was originally considered in \cite{Ben-Tal:2005:RFC:1246360.1246366}, describing the negotiation of flexible contracts between a retailer and a supplier in the presence of uncertain orders from customers. 
In particular, we analyze the performance of the approach proposed in this paper on a simplified version discussed in \cite{5986692} and in \cite{Vayanos2012459}.  
We remark that the considered numerical problem is such that the optimal solution of the multistage robust optimization problem can be assessed: this allows to evaluate the performance of the scenario with certificate approach.

The problem setting can be summarized as follows: a retailer received orders $\xi^t$ at the beginning of each time period $t\in\mathbb{T}=\left\{1,\dots,H-1\right\}$, $\xxi^{t}$ represents the demand history up to time $t$.
The demand needs to be satisfied from an inventory with filling level $s_{inv}^{t}$ by means of orders $x_{o}^{t}$ at a cost $d^t$ per unit of product.
Unsatisfied demand  may be backlogged at cost $p^t$ and inventory may be held in the warehouse with a unitary holding cost $h^t$. 
Lower and upper bounds on the orders $x_{o}^t$ ($\underline{x}_{o}^t$ and $\bar{x}_{o}^t$) at each period as well as on the cumulative orders $s_{co}^t$  ($\underline{s}_{co}^t$ and $\bar{s}_{co}^t$) up to period $t$ are imposed. We assume that there is no demand at time $t=1$ and that the demand at time $t$ lies within an interval centered around a nominal value $\bar{\xi}^{t}$ and uncertainty level $\rho\in\left[0,1\right]$ resulting in a box uncertainty set as follows: $\Xi = \times_{t\in \mathbb{T}} \left\{\xi^t\in\mathbb{R}: \left| \xi^t - \bar{\xi}^t\right|\leq \rho\bar{\xi}^t\right\}$. 
Denoting  with $x_{c}^t$ the retailer's cost at stage $t$, the problem can be modeled as follows
\begin{subequations}
\begin{eqnarray}
\textrm{RO}_H(\mathcal{COC}) {:=} & &\min_{x_{o}^t,x_{c}^{t},s_{co}^{t},s_{inv}^{t}}\left[ x_{c}^{1} +\max_{\xxi \in \Xi}  \sum_{t \in \mathbb{T}}x_{c}^{t+1}(\xxi^{t}) \right]\label{objese} \\
& & \textrm{s.t. } x_{c}^{1}  \geq  d^1 x_{o}^{1} + \max\left\{h^{1} s_{inv}^{1},-p^{1} s_{inv}^{1}\right\}  \label{conses10}\\
& &\qquad  x_{c}^{t+1}(\xxi^{t})  \geq  d^{t+1} x_{o}^{t+1}(\xxi^{t}) + \nonumber\\
& &\qquad  +\max\left\{h^{t+1} s_{inv}^{t+1}(\xxi^{t}),-p^{t+1} s_{inv}^{t+1}(\xxi^{t})\right\},\   t=1,\dots,H\!-\!2 \label{conses1}\\
& &\qquad x_{c}^{H}(\xxi^{H-1})  \geq    \max\left\{h^H s_{inv}^{H}(\xxi^{H-1}),-p^H s_{inv}^{H}(\xxi^{H-1})\right\}
\label{conses2} \\
& &\qquad s_{inv}^{2}(\xxi^{1})  =   s_{inv}^1 +  x_{o}^1 - \xi^{1}
\label{conses3bis}\\
& &\qquad s_{inv}^{t+1}(\xxi^{t})  =   s_{inv}^{t}(\xxi^{t-1}) +  x_{o}^{t}(\xxi^{t-1}) - \xi^{t}\ , \quad t=2,\dots, H-1
\label{conses3}\\
& &\qquad s_{co}^{2}(\xxi^{1})  =   s_{co}^{1} +  x_{o}^{1}
\label{conses4bis}\\
& &\qquad s_{co}^{t+1}(\xxi^{t})  =   s_{co}^{t}(\xxi^{t-1}) +  x_{o}^{t}(\xxi^{t-1})\ , \quad t=2,\dots, H-1
\label{conses4}\\
& &\qquad \underline{x}_{o}^{1}   \leq   x_{o}^{1} \leq \bar{x}_{o}^{1}, \quad
 \underline{s}_{co}^{1}   \leq  s_{co}^{1} \leq \bar{s}_{co}^{1}
\label{conses6bis}\\
& &\qquad \underline{x}_{o}^{t}   \leq   x_{o}^{t}(\xxi^{t-1}) \leq \bar{x}_{o}^{t}, \quad 
\underline{s}_{co}^{t}   \leq  s_{co}^{t}(\xxi^{t-1}) \leq \bar{s}_{co}^{t}, \quad t=2,\dots, H.
\label{conses6}
\end{eqnarray}
 \end{subequations} 
The  objective  function (\ref{objese}) corresponds to  minimize the worst-case cumulative cost. Constraints (\ref{conses10})-(\ref{conses1})-(\ref{conses2}) define the stage-wise {costs  $x_{c}^{t}(\xxi^{t})$, $t=1,\dots,H$} while constraints (\ref{conses3bis})-(\ref{conses3}) and (\ref{conses4bis})-(\ref{conses4}) respectively define the dynamics of the inventory level  and cumulative orders. Finally, constraints (\ref{conses6bis})-(\ref{conses6}) denote the lower and upper bounds on the instantaneous and cumulative orders. Notice that the decision process is  nonanticipative.


The corresponding multistage scenario with certificates  formulation, based on the extraction of $N$ samples of the uncertainty, becomes as follows
\begin{eqnarray*}
\textrm{SwC}_H^N(\mathcal{COC}) & {:=} & \min_{x_{o,i}^t,x_{c,i}^{t},s_{co,i}^{t},s_{inv,i}^{t}}  \gamma  \label{objesers} \\
& & \textrm{s.t. } \gamma \geq \left[ x_{c}^{1} + \sum_{t \in \mathbb{T}}x_{c,i}^{t+1}({\xxi^{t}}^{(i)}) \right]\ ,\quad  i=1,\dots,N  \label{epigraphss}\nonumber\\
& & x_{c}^{1}  \geq  d^1 x_{o}^{1} + \max\left\{h^{1} s_{inv}^{1},-p^{1} s_{inv}^{1}\right\}  \label{conses10as}\nonumber\\
& & x_{c,i}^{t+1}({\xxi^{t}}^{(i)})  \geq  d^{t+1} x_{o,i}^{t+1}({\xxi^{t}}^{(i)}) + \nonumber\\
& & +\max\left\{h^{t+1} s_{inv,i}^{t+1}({\xxi^{t}}^{(i)}),-p^{t+1} s_{inv,i}^{t+1}({\xxi^{t}}^{(i)})\right\},\   t=1,\dots,H\!-\!2,\ i=1,\dots,N  \label{conses1as}\nonumber\\
& & x_{c,i}^{H}({\xxi^{H-1}}^{(i)})  \geq    \max\left\{h^H s_{inv,i}^{H}({\xxi^{H-1}}^{(i)}),-p^H s_{inv,i}^{H}({\xxi^{H-1}}^{(i)})\right\},\  i=1,\dots,N   
\nonumber\label{conses2as} \\
& & s_{inv,i}^{2}({\xxi^{1}}^{(i)})  =   s_{inv}^1 +  x_{o}^1 - \xi^{1},\   i=1,\dots,N\nonumber
\label{conses3bisas}\\
& & s_{inv,i}^{t+1}({\xxi^{t}}^{(i)}) 	\! =\!   s_{inv,i}^{t}({\xxi^{t-1}}^{(i)})\! +\!  x_{o,i}^{t}({\xxi^{t-1}}^{(i)})\! -\! {\xi^{t}}^{(i)}, \ t=\!2,\!\dots\!,\!H-1,\ i=\!1,\!\dots\!,\!N\nonumber
\label{conses3as}\\
& & s_{co,i}^{2}({\xxi^{1}}^{(i)})  =   s_{co}^{1} +  x_{o}^{1},\  i=1,\dots,N \nonumber
\label{conses4bisas}\\
& & s_{co,i}^{t+1}({\xxi^{t}}^{(i)}) \! = \!  s_{co,i}^{t}({\xxi^{t-1}}^{(i)})\! +\!  x_{o,i}^{t}({\xxi^{t-1}}^{(i)}),\  t=2,\!\dots,\!H-1,\ i=1,\!\dots,\!N 
\nonumber\label{conses4as}\\
& & \underline{x}_{o}^{1}   \leq   x_{o}^{1} \leq \bar{x}_{o}^{1},\quad
\underline{s}_{co}^{1}   \leq  s_{co}^{1} \leq \bar{s}_{co}^{1}
\nonumber\label{conses6bisas}\\
& & \underline{x}_{o}^{t}   \leq   x_{o,i}^{t}({\xxi^{t-1}}^{(i)}) \leq \bar{x}_{o}^{t}, \quad t=2,\dots, H,\ i=1,\dots,N 
\nonumber\label{conses5as}\\
& & \underline{s}_{co}^{t}   \leq  s_{co,i}^{t}({\xxi^{t-1}}^{(i)}) \leq \bar{s}_{co}^{t}, \quad t=2,\dots, H,\  i=1,\dots,N.\nonumber
\label{conses6as}
\end{eqnarray*}

We consider specific instances of  problem $\textrm{RO}_H(\mathcal{COC})$ as  summarized in Table \ref{data} under the assumption of  two-stage ($H=2$)  and a five-stage ($H=5$) time horizons.
The data presents some slight modifications of the data presented in \cite{Vayanos2012459}. 

\begin{table}[ht!]
	\centering
		\begin{tabular}{ll}
		\hline
		Parameters 		 &   $\textrm{RO}_H(\mathcal{COC})$  \\
			\hline
$H$ & 2 $/$ 5 \\
$(p^t,d^t,h^t)$ & (11,1,10)\\
$(s_{inv}^{1})$ & 0 \\
($\underline{x}_{o}^{t},\bar{x}_{o}^{t}$) & $(0,\infty)$ \\
($\underline{s}_{co}^{1},\dots,\underline{s}_{co}^{H-1}$) & $(47,134,188,429)$ \\
($\bar{s}_{co}^{1},\dots,\bar{s}_{co}^{H-1}$) & $(94,248,370,586)$ \\
$\bar{\xi}^t$, $t=1,\dots,H-1$ & $100 \left(1 + \frac{1}{2}\sin\left(\frac{\pi\left(t- 2\right)}{6}\right)\right)=(75,100,125,143.3013)$ \\
$\rho $ & $30\%$ \\
\hline
		\end{tabular}
		\caption{Input data for the inventory management problem.}
		\label{data}
\end{table}

We define optimality gaps of the problem $\textrm{SwC}_H^N(\mathcal{COC})$ as
{
\begin{equation}
optimality\  gap :=\frac{\inf \textrm{SwC}_H^N(\mathcal{COC}) - \inf \textrm{RO}_H(\mathcal{COC}) }{\inf \textrm{RO}_H(\mathcal{COC})}.
\label{optgap}
\end{equation}
}
We note that the optimality gap in (\ref{optgap}) can be computed, since problem $\textrm{RO}_H(\mathcal{COC})$ can be solved exactly by using a scenario tree that consists of the vertices of $\Xi$ reported in Tables \ref{vertices2} and~\ref{vertices5}. 

\begin{table}[ht!]
	\centering
		\begin{tabular}{ll}
		\hline
		 		 &   $t=2$    \\
			\hline
1 &		52.5 \\
2	&  97.5 \\
\hline
		\end{tabular}
		\caption{Vertices of $\Xi$ for the management inventory problem in the two-stage case ($H=2$).}
		\label{vertices2}
\end{table}

\begin{table}[ht!]
	\centering
		\begin{tabular}{lllll}
		\hline
		 		 &  $t=1$ & $t=2$ & $t=3$ & $t=4$    \\
			\hline
1 &		52.5	& 70	& 87.5 &	100 \\
2	&  52.5 &	70	& 87.5 &	186 \\
3	&  52.5 &	70 &	163 &	100 \\
4 &		52.5 &	70 &	163 &	186 \\
5	& 	52.5 &	130 &	87.5 &	100 \\
6 &		52.5 &	130 &	87.5 &	186 \\
7 &		52.5 &	130 &	163 &	100 \\
8 &		52.5 &	130 &	163 &	186 \\
9 &		97.5 &	70 &	87.5 &	100 \\
10 &		97.5 &	70 &	87.5 &	186 \\
11 &		 97.5 &	70	& 163 &	100 \\
12 &	 97.5 &	70	 & 163 &	186 \\
13	&  97.5	& 130	& 87.5 &	100 \\
14 &	 97.5 &	130	& 87.5 &	186 \\
15 &	 97.5 &	130	& 163	& 100 \\
16 &		97.5 &	130 &	163	& 186 \\
\hline
		\end{tabular}
		\caption{Vertices of $\Xi$ for the management inventory problem in the five-stage case ($H=5$).}
		\label{vertices5}
\end{table}

To assess the performance of our approach, we compute the \textit{empirical violation probability} of a given solution $(x^1,\gamma)$, defined as follows:

\begin{displaymath}
\widehat{\textrm{V}}_H(x^1,\gamma)\!\! :=\!\! \frac{1}{L}\!\!\sum_{\ell=1}^{L}\!\! \left\{
\begin{array}{l}
 \mbox{\boldmath $\hat{\xi}$}^{H-1}_{\ell}\in \Xi \textrm{ s.t.} \nonumber\\
\left\{ \begin{array}{l}
   \frac{1}{N}\sum_{i=1}^{N}\mathds{1}[\not\exists\  x^{t}(\mbox{\boldmath $\hat{\xi}$}^{t-1(i)}_{\ell}) \in\mathbb{R}^{n_t}_+,\ t=2,\dots,H \textrm{ satisfying }\nonumber \\
  {c^1}^{\top} x^1+{c^{2}}^\top\!\!\!\left(\mbox{\boldmath $\hat{\xi}$}^{1(i)}_{\ell}\right) x^{2}\!\left(\mbox{\boldmath $\hat{\xi}$}^{1(i)}_{\ell}\right)+\cdots+
{c^{H}}^\top\!\!\!\left(\mbox{\boldmath $\hat{\xi}$}^{H-1(i)}_{\ell}\right) x^{H}\!\left(\mbox{\boldmath $\hat{\xi}$}^{H-1(i)}_{\ell}\right) \le \gamma\\
   T^{t-1}(\mbox{\boldmath $\hat{\xi}$}^{t-1(i)}_{\ell})x^{t-1}(\mbox{\boldmath $\hat{\xi}$}^{t-2(i)}_{\ell}) + W^t x^{t}(\mbox{\boldmath $\hat{\xi}$}^{t-1(i)}_{\ell})= h^t(\mbox{\boldmath $\hat{\xi}$}^{t-1(i)}_{\ell}),\  t=2,\dots,H]  \nonumber 
\end{array} \right. 
\end{array}  \right\},
\label{empviol}
\end{displaymath}
where $\mathds{1}$ is the indicator function counting the number of scenarios where the constraints are not satisfied and $\left\{\mbox{\boldmath $\hat{\xi}$}^{H-1(i)}_{\ell}\right\}_{\ell=1,\dots,L}$ is a sequence of  $L$ independent sets distributed according to $\pr$. The ${\ell}-th$-sample  is composed by  scenarios $\mbox{\boldmath $\hat{\xi}$}^{H-1(1)}_{\ell},\dots,\mbox{\boldmath $\hat{\xi}$}^{H-1(N)}_{\ell}$.  Notice that these samples are independent of the  $N$   samples  $\xxi^{H-1(i)}$ used in problem (\ref{MSwC_N}) to obtain solution $(\hat{x}^1,\hat{\gamma})$.
\vskip 4mm
\noindent
The numerical results are obtained as follows: 
\begin{itemize}
\item[-] we fix a confidence level of $\beta=0.1\%$ for the constraint sampling;   
	\item[-] we select the target violation probability $\epsilon=0.00025,0.0005,0.001,0.005,0.01,0.05,0.1,0.2,0.3,$ and compute the corresponding sample size $N=N(\epsilon,\beta)$ according to formula (\ref{Al})  rounded up to the next integer;
\item[-] we solve 100 instances of problem $SwC_N$ each based on a different number $N$ of independent samples as computed in the previous point;
\item[-] for each instance, we compute the optimality gap given in formula (\ref{optgap}) and empirical violation probability given in formula (\ref{empviol}) with $L=1000$;
\item[-]  we compute statistics over the 100 istances.
\end{itemize}

The problems derived from the case study have been formulated and solved under \textit{AMPL} environment along \textit{CPLEX} {12.5.1.0} solver.
All  computations have been performed on a $64$-bit machine with 12 GB of RAM and a Intel Core i7-3520M CPU $2.90$ GHz processor.

\begin{figure}[ht!]
\centering
\includegraphics[width=\textwidth]{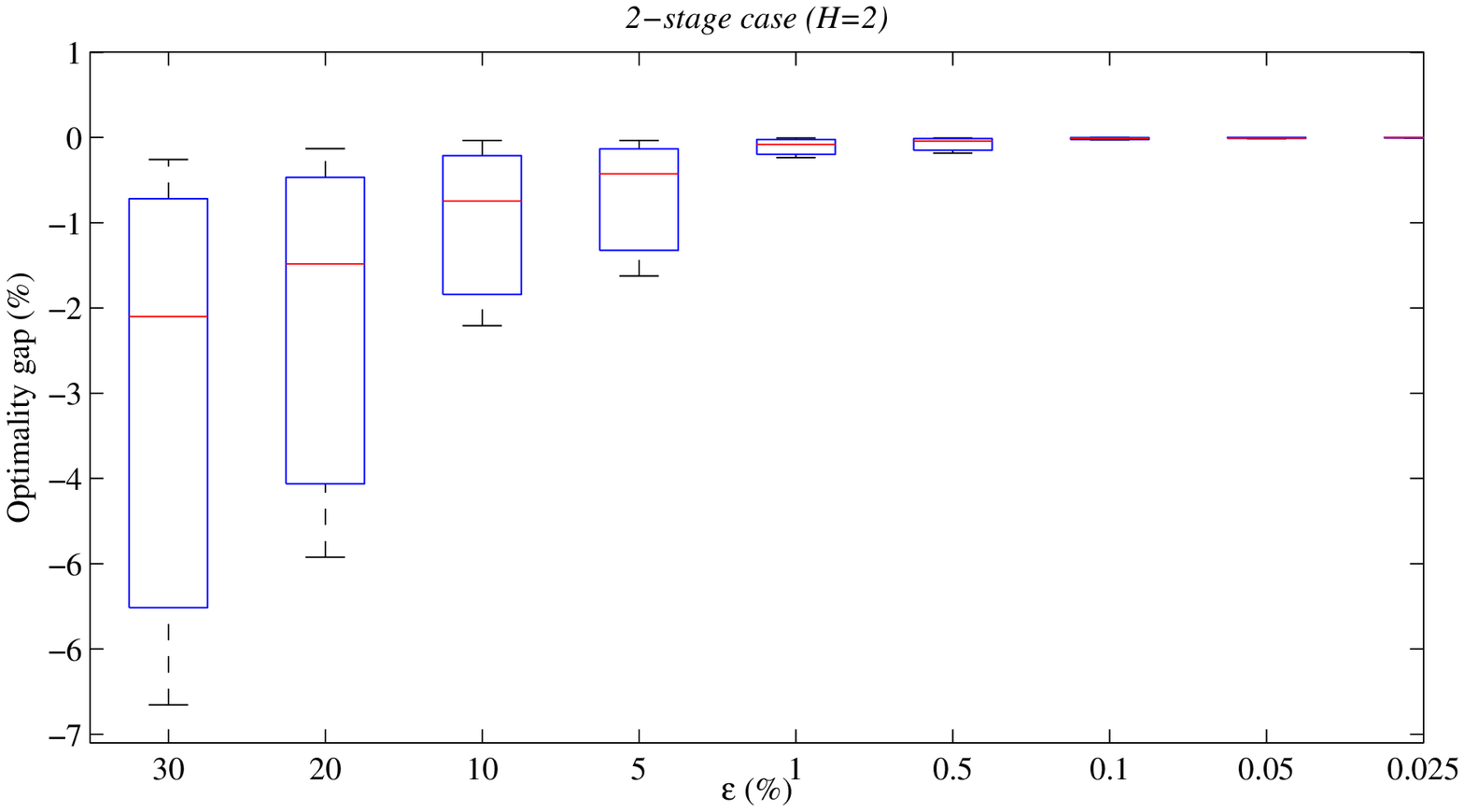}
\vskip3truemm\noindent
\caption{Optimality gaps for $\textrm{SwC}_2^N(\mathcal{COC})$ (boxes and whiskers) for decreasing values of $\epsilon$  for the two-stage ($H=2$) case.}
\label{optimalitygap2}
\end{figure}

\begin{figure}[ht!]
\centering
\includegraphics[width=\textwidth]{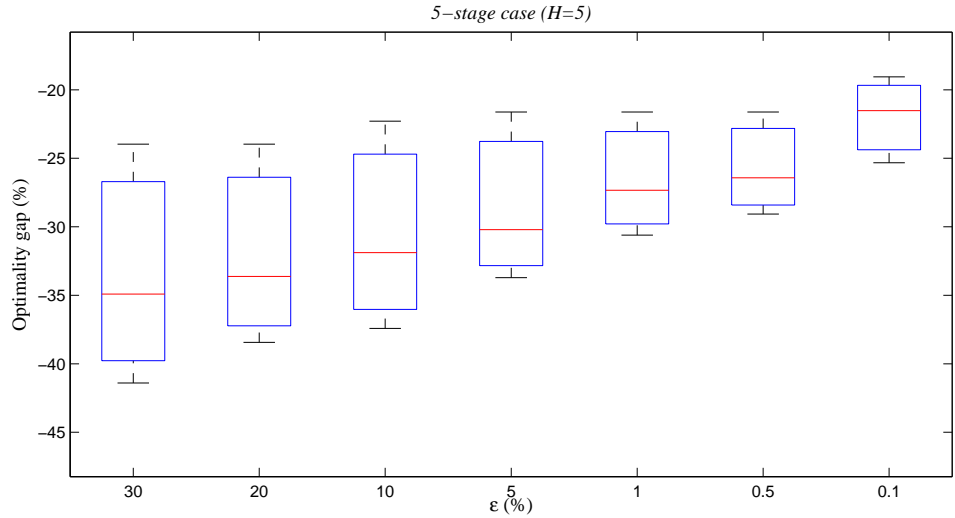}
\vskip3truemm\noindent
\caption{Optimality gaps for $\textrm{SwC}_5^N(\mathcal{COC})$ (boxes and whiskers) for decreasing values of $\epsilon$ for the five-stage ($H=5$)  case.}
\label{optimalitygap5nonanticip}
\end{figure}

First, we evaluate the performance of the scenario with certificates {$\textrm{SwC}_H^N(\mathcal{COC})$} approach.
Figures \ref{optimalitygap2} and \ref{optimalitygap5nonanticip} display the optimality gaps of problems $\textrm{RO}_H(\mathcal{COC})$ with respect to $\inf \textrm{RO}_H(\mathcal{COC})$ for different values of violation probability $\epsilon$ ($\%$) ranging from $30\%$ up to $0.025\%$ respectively for the two ($H=2$) and five-stage  ($H=5$) cases. {Number of samples $N$, constraints and variables of the corresponding optimization models  are reported in Table \ref{sampleN}.}

\begin{table}[ht!]
	\centering
		\begin{tabular}{llllll}
		\hline
		  $\epsilon$ ($\%$) & $N$ & $\#$ of const. ($H=2$) & $\#$ of var. ($H=2$) & $\#$ of const. ($H=5$) & $\#$ of var. ($H=5$) \\
			\hline
30 & 63& 756 & 442 & 2772 & 1198 \\
20 & 95 & 1140 & 666 & 4180 & 1806 \\
10 & 189 & 2268 & 1324 & 8316 & 3592 \\
5 & 377 & 4524 & 2640 & 16588 & 7164 \\
1 & 1884 & 22608 & 13189 & 82896 & 35797 \\
0.5 & 3768 & 45216  & 26377 & 165792 & 71593\\
0.1 & 18838 & 226056 & 131867 & 414433 & 339085 \\
0.05 & 37676  & 452112 & 263733 & 828869 & 678169 \\ 
0.025 & 75352 & 904224  & 527465 & 1657741 &  1356337 \\
\hline
		\end{tabular}
		\caption{Number of samples $N$, constraints and variables for decreasing values of  $\epsilon$ ($\%$) both in the two-stage {(columns 3 and 4)} and five stage {(columns 5 and 6)} cases.}
		\label{sampleN}
\end{table}

From the results shown in Figure \ref{optimalitygap2} and \ref{optimalitygap5nonanticip} we can observe that the variance of {$\textrm{SwC}_H^N(\mathcal{COC})$} decreases substantially as $\epsilon $ decreases as well as the optimality gaps passing from $-2\%$ (in average) to $-10^{-3}\%$ for the two-stage case and from  $-34\%$ (in average) to $-21\%$ for the five-stage case.
It should be observed that, for the same given level of allowed violation $\epsilon$, the $\textrm{SwC}_5^N(\mathcal{COC})$ cost will always be lower than the solution returned by the approach in \cite{Vayanos2012459} (the reader is referred to the example proposed in that paper for comparison). This is due to the conservatism  introduced in \cite{Vayanos2012459} by the fact that the solution is constrained to variables for stages $2,\ldots,5$ with a special structure, and it is the reason why we have larger optimality gaps. We stress that this is a desirable feature, since we find a better result using the same level of probability. 

\begin{figure}[ht!]
\centering
\includegraphics[width=\textwidth]{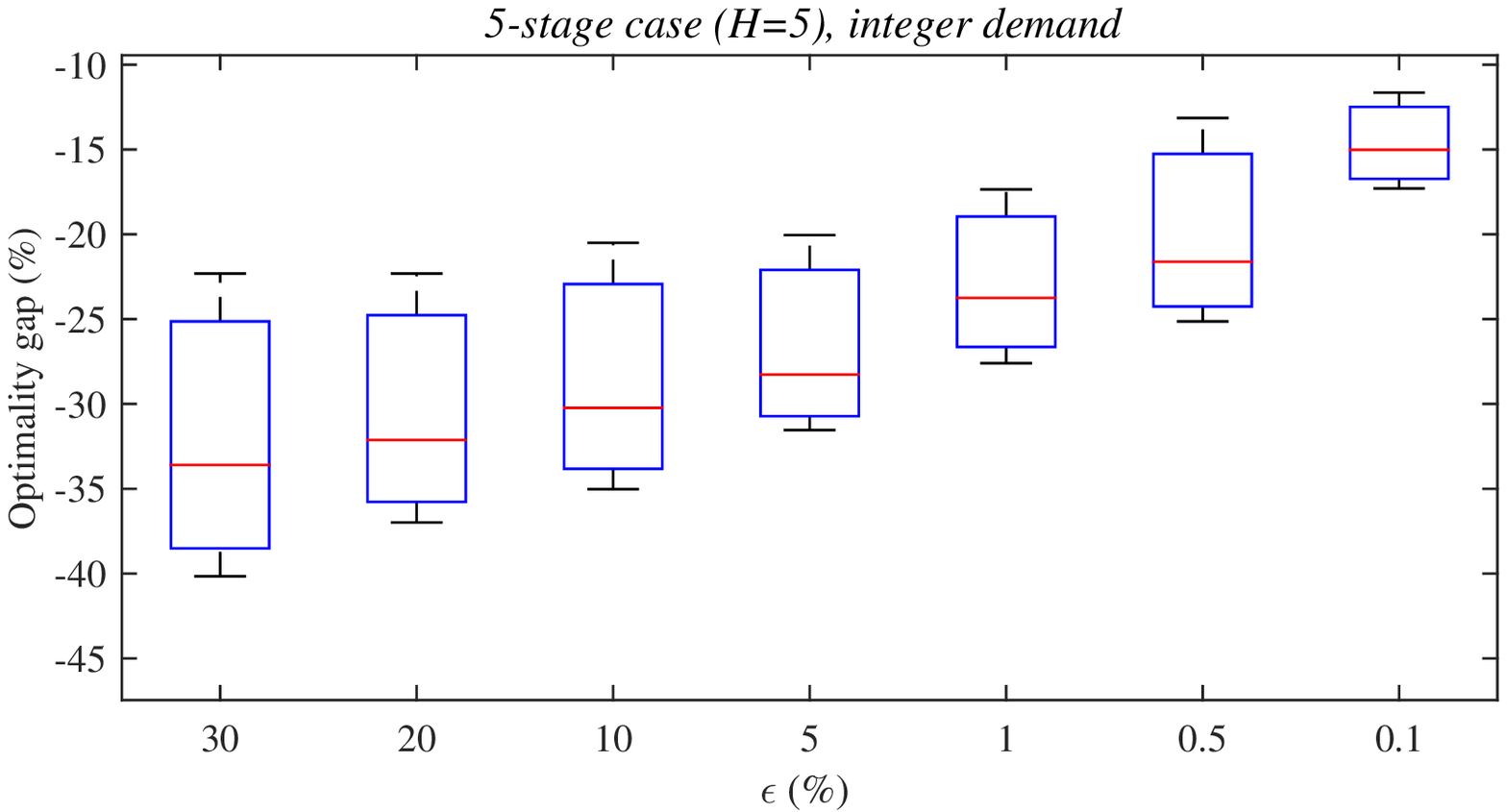}
\vskip3truemm\noindent
\caption{Optimality gaps for $\textrm{SwC}_5^N(\mathcal{COC})$ (boxes and whiskers) for decreasing values of $\epsilon$ for the five-stage ($H=5$)  case with integer demand over stages.}
\label{optimalitygap5nonanticipINTEGER}
\end{figure}

We also note that, since the uncertainty lies in continuos intervals, we have a  probability close to zero to get twice the same sample. Consequently,  the nonanticipativity constraints in problem 
$\textrm{SwC}_5^N(\mathcal{COC})$ are not required to be satisfied by our data, with resulting lower costs in term of objective function values. This was not the case with $H=2$.
Based on this observation, we performed a second computational test, shown in Figure \ref{optimalitygap5nonanticipINTEGER}, where the demand $\xi^t$ $t=1,\dots,4$ is assumed to take only integer values in the intervals $[52.5, 97.5]$,	$[70,	130]$, $[87.5,163]$, $[100,186]$. In this way, we increase the probability of having repeated samples, thus  enforcing nonaticipativity constraints. Results show that the optimality gaps are now reduced passing from  $-33\%$ (on average) to $-15\%$ for the five-stage case.

\begin{figure}[ht!]
\centering
\includegraphics[width=\textwidth]{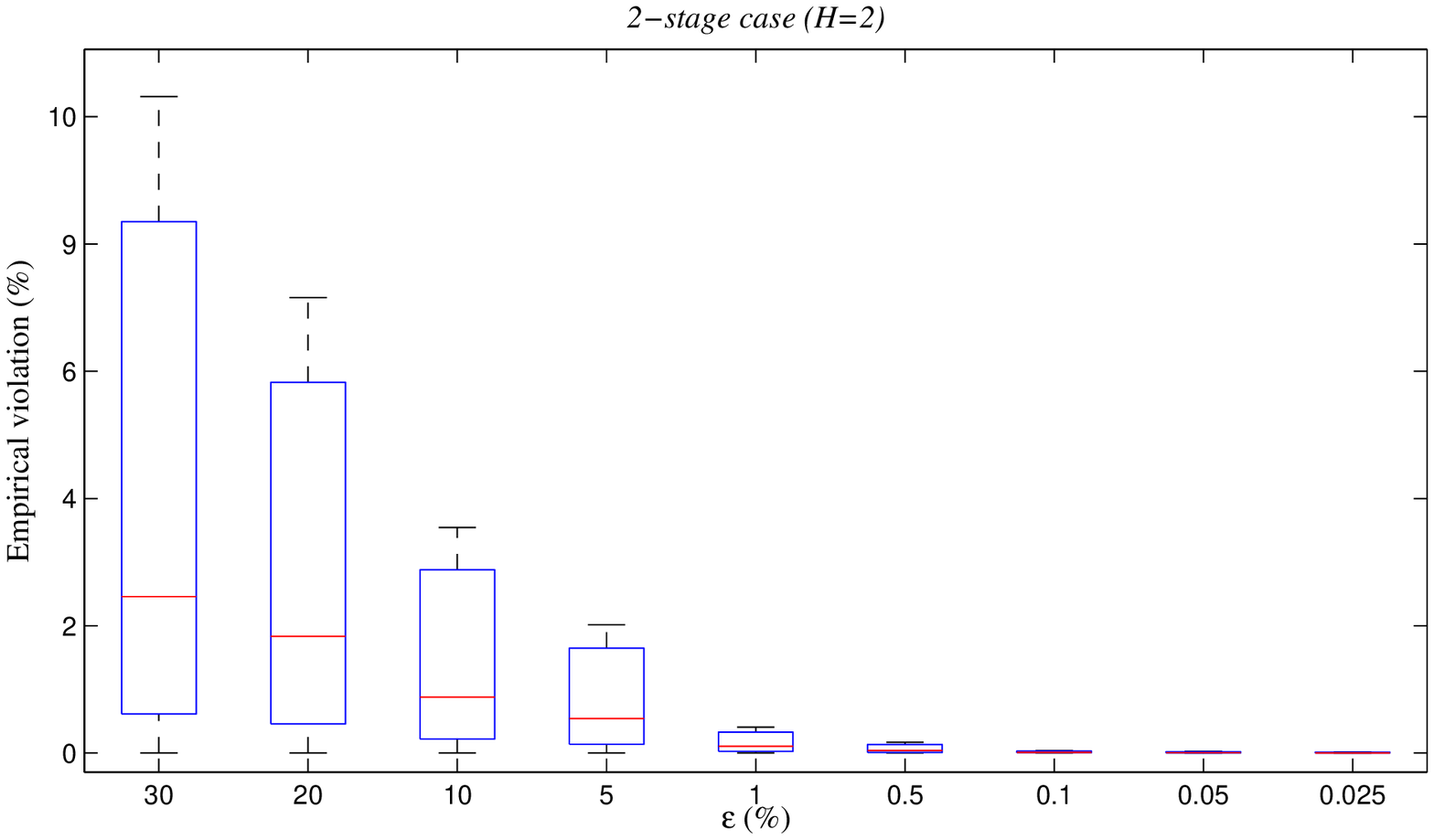}
\vskip3truemm\noindent
\caption{Empirical violation probability  for $\textrm{SwC}_2^N(\mathcal{COC})$ (boxes and whiskers) for increasing values of $\epsilon$ for the two-stage ($H=2$)  case.}
\label{violprob2}
\end{figure}

\begin{figure}[ht!]
\centering
\includegraphics[width=\textwidth]{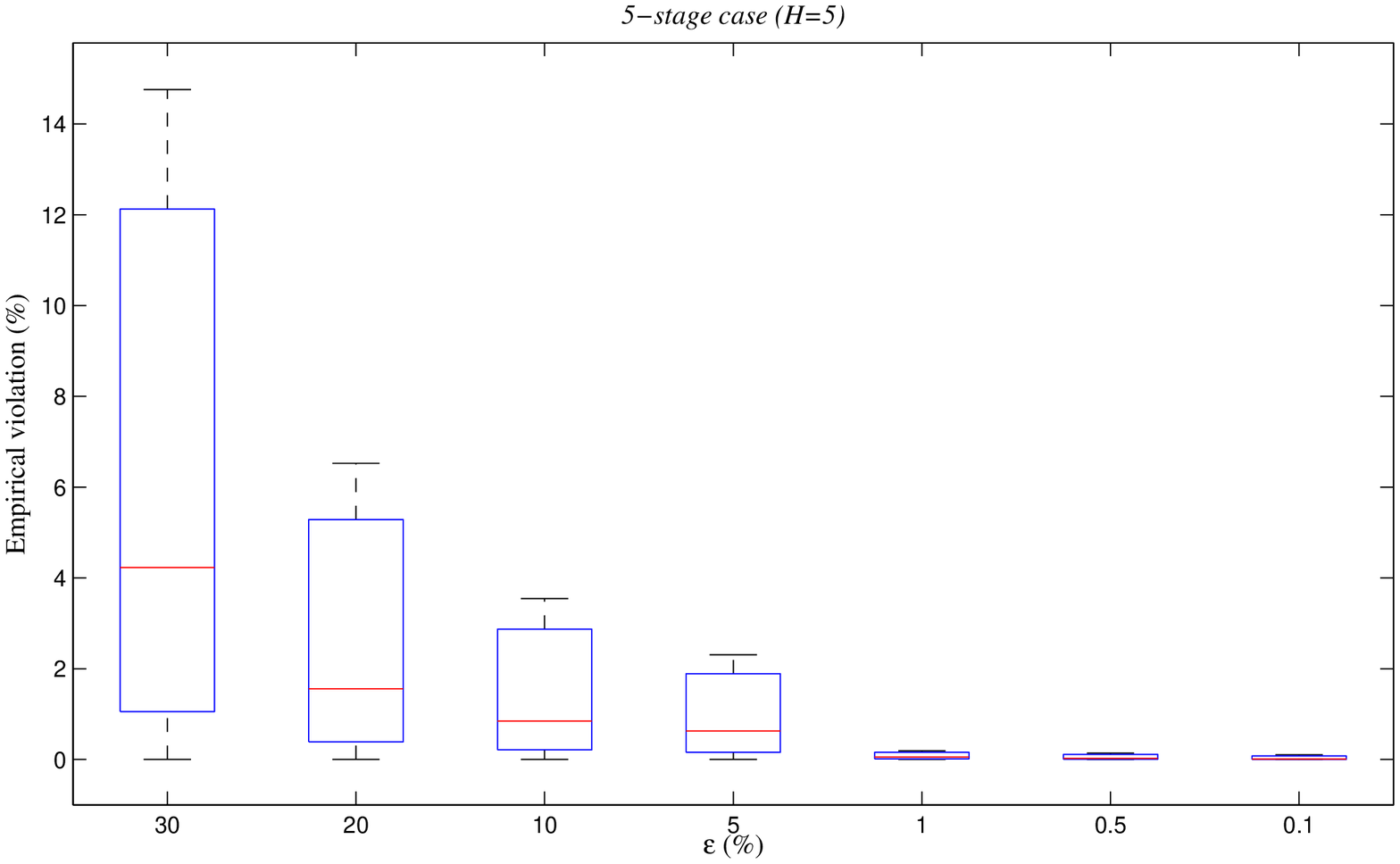}
\vskip3truemm\noindent
\caption{Empirical violation probability  for $\textrm{SwC}_5^N(\mathcal{COC})$ (boxes and whiskers) for decreasing values of $\epsilon$ for the five-stage ($H=5$)  case.}
\label{violprob5}
\end{figure}

In Figures \ref{violprob2} and \ref{violprob5}, we plot the distribution of the empirical violation probability as function of  $\epsilon$, for the two-stage ($H=2$) case and the five-stage ($H=5$) case.   As expected, as $\epsilon$ decreases, the violation converges to 0. We also note that the empirical violation probability is smaller than $\epsilon$ in all the considered cases. 

\begin{figure}[ht!]
\centering
\includegraphics[width=\textwidth]{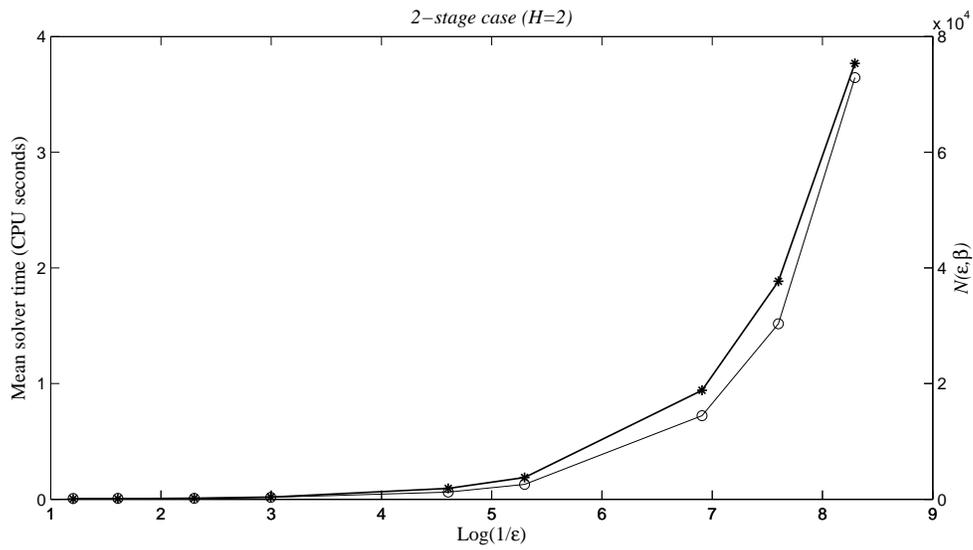}
\vskip3truemm\noindent
\caption{Mean solver times (solid lines) and number of samples (dashed lines) as a function of $\ln(1/\epsilon)$ for problem $\textrm{SwC}_2^N(\mathcal{COC})$  for the two-stage ($H=2$) case.}
\label{time2}
\end{figure}

\begin{figure}[ht!]
\centering
\includegraphics[width=\textwidth]{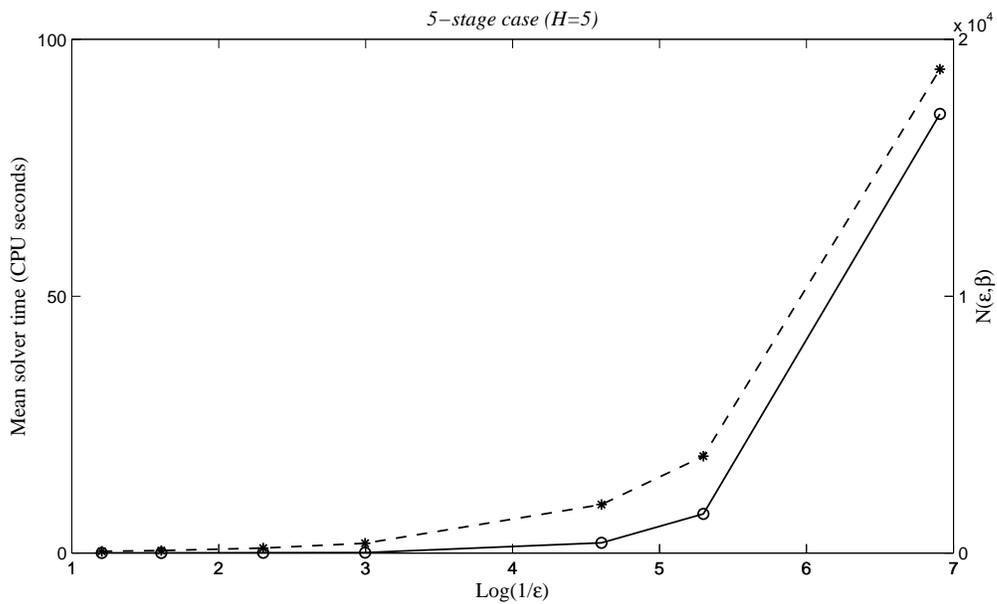}
\vskip3truemm\noindent
\caption{Mean solver times (solid lines) and number of samples (dashed lines) as a function of $\ln(1/\epsilon$) for problem $\textrm{SwC}_5^N(\mathcal{COC})$  for the five-stage ($H=5$)  case.}
\label{time5}
\end{figure}

Finally, Figures \ref{time2} and \ref{time5} show the average solver time (solid lines) and number of samples (dashed lines) for problems $\textrm{SwC}_2^N(\mathcal{COC})$ and $\textrm{SwC}_5^N(\mathcal{COC})$ as a function of $1/\epsilon$. 
Notice that the required number of samples obtained using formula  (\ref{Al}), corresponding to a prescribed level of violation probability does not depend on the number of stages and dimension of the certificates variables. Consequently, the number of samples shown in Table \ref{sampleN} are the same both for the two and five stage cases. In particular, they are considerably lower than those used in \cite{Vayanos2012459}, where the number of samples depends on the size of the basis and on the number of decision variables at each stage. On the other hand, we should remark that the number of variables used in our approach is larger, due to the introduction of sample-dependent certificates.

\subsection{Bounds for the Inventory management with cumulative orders}
In this section, we evaluate possible relaxations  to problem $\textrm{RO}_H(\mathcal{COC})$ as described in Section \ref{boundsMultistageRO}. In particular we consider the multistage wait-and-see problem $\textrm{RWS}_H(\mathcal{COC})$ for problem $\textrm{RO}_H(\mathcal{COC})$, and  the robust two-stage relaxation problem {$\textrm{RT}_H(\mathcal{COC})$} where the nonanticipativity constraints are relaxed in stages $2,\dots,H$.
Again, we remark that for the case at hand these two problems can be computed exactly by considering only the vertices of $\Xi$.
Similarly to formula (\ref{optgap}), we define  optimality gaps of the problem {$\textrm{RWS}_H(\mathcal{COC})$} as
\begin{equation}
(optimality\  gap)_{\textrm{RWS}_H(\mathcal{COC})} :=\frac{\inf \textrm{RWS}_H(\mathcal{COC}) - \inf \textrm{RO}_H(\mathcal{COC})}{\inf \textrm{RO}_H(\mathcal{COC})},
\label{optgapbounds}
\end{equation}
and similarly for $\textrm{RT}_H(\mathcal{COC})$. 
Both optimality gaps of  {$\textrm{RWS}_5(\mathcal{COC})$}  and of  $\textrm{RT}_5(\mathcal{COC})$ turned out to be equal to $ -0.170171593$,  passing from an objective function value of $2207.554108$ for $\textrm{RO}_5(\mathcal{COC})$ to $1831.891109$. Consequently the  Robust Value of Perfect Information {$\textrm{RVPI}_5$} is $375.66299$.

\begin{figure}[ht!]
\centering
\includegraphics[width=\textwidth]{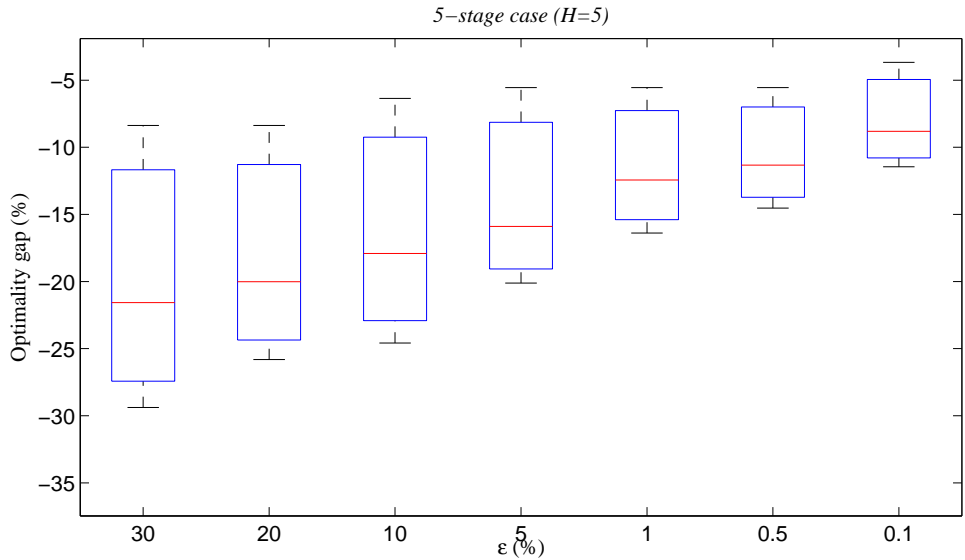}
\vskip3truemm\noindent
\caption{Optimality gaps for {$\textrm{SwCT}_5^N(\mathcal{COC})$} with respect to the robust two-stage
relaxation problem $\textrm{RT}_5^N(\mathcal{COC})$, (boxes and whiskers) for decreasing values of $\epsilon$ for the five-stage ($H=5$)  case.}
\label{optimalitygap5}
\end{figure}

\begin{figure}[ht!]
\centering
\includegraphics[width=\textwidth]{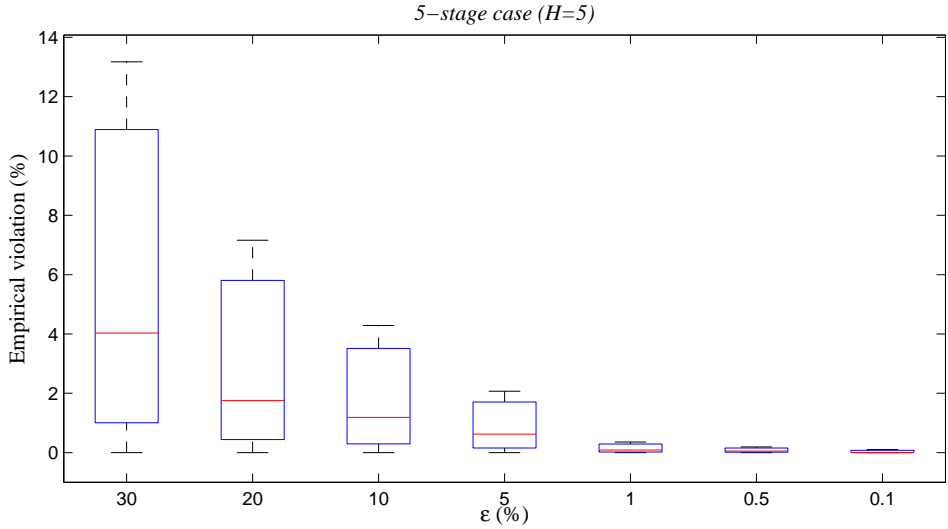}
\vskip3truemm\noindent
\caption{Empirical violation probability  for $\textrm{SwCT}_5^N(\mathcal{COC})$ (boxes and whiskers) with respect to the two-stage relaxation problem $\textrm{RT}_5^N(\mathcal{COC})$,  for decreasing values of $\epsilon$ in the five-stage ($H=5$)  case.}
\label{violprob5relax}
\end{figure}

Figure \ref{optimalitygap5} displays optimality gaps for the the two-stage relaxation scenario with
certificates  problem $\textrm{SwCT}_5^N(\mathcal{COC})$ with respect to the robust two-stage relaxation problem $\textrm{RT}_5^N(\mathcal{COC})$   for different values of violation probability $\epsilon$ ($\%$) ranging from $30\%$ up to $0.01\%$. From the results we can observe that the variance of $\textrm{SwCT}_5^N(\mathcal{COC})$ decreases substantially as $\epsilon $ decreases as well as the optimality gaps passing  from  $-21\%$ (in average) to $-8\%$. Note that the smaller values of optimality gaps compared to the ones obtained in Figure \ref{optimalitygap5nonanticip} for problem $\textrm{SwC}_5^N(\mathcal{COC})$ are mainly due to the fact that in problem $\textrm{SwCT}_5^N(\mathcal{COC})$ the nonanticipativity constraints are relaxed.
Finally Figure \ref{violprob5relax} refers to the empirical violation probability of $\textrm{SwCT}_5^N(\mathcal{COC})$ with respect to the robust two stage relaxation problem {$\textrm{RT}_5(\mathcal{COC})$}.  As expected as $\epsilon$ decreases it converges to 0. We again note that the empirical violation probability is smaller than $\epsilon$ in all the cases considered. 

\section{Conclusions}
\label{sec:Conclusions}

In this paper  probabilistic guarantees  for constraint sampling in multistage convex robust  optimization problems have been proposed.  
The scenario with certificates approach has been considered to treat   the dynamic nature of  convex multistage robust optimization problems.   
A multistage robust convex optimization problem  has been proved to be  equivalent to a convex robust optimization problem with certificates and  a bound on the probability of violation  for the scenario  with certificates approach  has been provided.   
The proposed approach has the important advantage to avoid the conservative use of parametrization through decision rules proposed in literature,  implying a strong reduction of the number of samples required to satisfy the required level of reliability. This is due to the fact that the required number of samples does not depend on the number of stages and dimension of the certificates variables. 
Numerical results on a case study taken from the literature show the efficiency of the proposed approach.

\appendix
\section{Proof of Theorem \ref{the2stagesamplecomplexity}}
 
We first prove the convexity of the set $\mathcal{X}_{\textrm{RwC}_2}(\xi^1)$ defined above.
Given $\hat{x}^1,\tilde{x}^1\in\mathcal{X}_{\textrm{RwC}_2}(\xi^1)$, then there exist $\hat{x}^2(\xi^1),\tilde{x}^2(\xi^1)$ such that
\begin{displaymath}
 \left\{
\begin{array}{l}
 T^1(\xi^1)\hat{x}^1 + W^2(\xi^{1}) \hat{x}^2(\xi^1) =  h^2(\xi^1),\  c^{{1}^{\top}}{\hat{x}^1} + c^{{2}^{\top}}(\xi^1) \hat{x}^2(\xi^1)\leq \gamma   \\
\label{convexity}
 T^1(\xi^1)\tilde{x}^1 + W^2(\xi^{1}) \tilde{x}^2(\xi^1)  =  h^2(\xi^1),\  c^{{1}^{\top}}{\tilde{x}^1} + c^{{2}^{\top}}(\xi^1) \tilde{x}^2(\xi^1)\leq \gamma\ .
\end{array}
\right.
\end{displaymath}
Consider now $x^{1\lambda}:=\lambda \hat{x}^1 + (1 - \lambda)\tilde{x}^1,$ with $\lambda\in[0,1]$, and let $x^{2\lambda}:=\lambda \hat{x}^2(\xi^1) + (1 -\lambda)\tilde{x}^2(\xi^1)$. 
Then
\begin{eqnarray}
 T^1(\xi^1)x^{1\lambda} + W^2(\xi^{1}) x^{2\lambda} & = &  T^1(\xi^1)\left(\lambda \hat{x}^1 + (1 - \lambda)\tilde{x}^1\right) + W^2(\xi^{1}) \left(\lambda \hat{x}^2(\xi^1) + (1 -\lambda)\tilde{x}^2(\xi^1)\right) \nonumber\\
&=&\lambda\left(T^1(\xi^1)\hat{x}^1 + W^2(\xi^{1}) \hat{x}^2(\xi^1) \right)  +  (1 -\lambda)\left( T^1(\xi^1)\tilde{x}^1 + W^2(\xi^{1}) \tilde{x}^2(\xi^1)\right) \nonumber\\
& = & \lambda h^2(\xi^1) + (1 -\lambda) h^2(\xi^1) \nonumber\\
&= & h^2(\xi^1),\nonumber 
\end{eqnarray}
and
\begin{eqnarray}
c^{{1}^{\top}}x^{1\lambda} + c^{{2}^{\top}}(\xi^1)  x^{2\lambda} & = & c^{{1}^{\top}}\left(\lambda \hat{x}^1 + (1 - \lambda)\tilde{x}^1\right) + c^{{2}^{\top}}(\xi^1)\left(\lambda \hat{x}^2(\xi^1) + (1 -\lambda)\tilde{x}^2(\xi^1)\right)\nonumber \\
& = &  \lambda\left(c^{{1}^{\top}}{\hat{x}^1} + c^{{2}^{\top}}(\xi^1) \hat{x}^2(\xi^1) \right) + (1 -\lambda)\left(c^{{1}^{\top}}{\tilde{x}^1} + c^{{2}^{\top}}(\xi^1) \tilde{x}^2(\xi^1)\right)\nonumber\\
& \leq &  \lambda \gamma + (1 -\lambda)\gamma \nonumber\\
& = & \gamma,\nonumber
\end{eqnarray}
which proves the convexity of $\mathcal{X}_{\textrm{RwC}_2}(\xi^1)$.
From Theorem \ref{lemma1} we observe that the condition $({x^1},\gamma)\in\mathcal{X}_{\textrm{RwC}_2}(\xi^1)$ is equivalent to $({x^1},\gamma)\in\mathcal{X}_{\textrm{RO}_2}(\xi^1)$ and that the problem $\textrm{RwC}_2$ is equivalent to $\textrm{RO}_2$ given by
\begin{eqnarray}
\label{ro_twoquadris}
\min_{ {x^1},\gamma} &&  \gamma  \\
& \textrm{s.t. }&			 A {x^1}=h^1,\quad {x^1}\geq 0\nonumber\\
& &c^{{1}^{\top}}{x^1} + \min_{{x^2}(\xi^1)} c^{{2}^{\top}}\left(\xi^{1}\right){x^2}(\xi^1)\leq  \gamma \nonumber \\
& &   T^1(\xi^1){x^1} + W^2(\xi^{1}) {x^2}(\xi^1)= h^2(\xi^1),\quad \quad {x^2}(\xi^1)\geq 0,\quad	\forall \xi^1\in\Xi^1. \nonumber
\end{eqnarray}
For the convexity of $\mathcal{X}_{\textrm{RwC}_2}(\xi^1)$ it follows that the following functions computed at the optimum of ${x^2}(\xi^1)$, say ${x^2}^{\ast}(\xi^1)$, $$ T^1(\xi^1){x^1} + W^2(\xi^{1}) {x^2}^{\ast}(\xi^1)= h^2(\xi^1),$$ and $$c^{{1}^{\top}}{x^1} +  c^{{2}^{\top}}\left(\xi^{1}\right){x^2}^{\ast}(\xi^1)\leq  \gamma$$ are convex in $ {x^1}$ for given $\xi^1$. Hence, the problem (\ref{ro_twoquadris}) is a robust convex optimization problem. 
Then, we construct its scenario counterpart
\begin{eqnarray}
\label{ro_two5}
\min_{ {x^1},\gamma}&&  \gamma  \\
& \textrm{s.t. }&			  A {x^1}=h^1,\quad {x^1}\geq 0\nonumber\\
& & c^{{1}^{\top}}{x^1} + \min_{{x^2_i}} c^{{2}^{\top}}\left({\xi^{1}}^{(i)}\right){x^2_i}\leq  \gamma\nonumber \\
& &  T^1({\xi^{1}}^{(i)}){x^1} + W^2({\xi^{1}}^{(i)}) {x^2_i}= h^2({\xi^{1}}^{(i)}),\quad \quad {x^2_i}\geq 0,\quad	i=1,\dots,N,\nonumber
\end{eqnarray}
where the subscript $i$ for the variables ${x^2_i}$ highlights that the different minimization problems are independent. Finally,  we note that the  problem (\ref{ro_two5}) corresponds to the  problem $\textrm{SwC}^N_2$ and the thesis follows from \cite{doi:10.1137/07069821X}. \qed

\section{Proof of Theorem \ref{theomultistageequiv}}

We first note that Problem {RO$_H$} can be rewritten in epigraph form, by introducing the additional variable $\gamma$, as follows
\begin{eqnarray}
\label{ROH_bbb}
\textrm{RO}_H& = & \! \min_{{x^1} ,\gamma}  \gamma \\
& \textrm{s.t. }&A{x^1} =h^1,\quad {x^1} \geq 0\nonumber\\
& &c^{{1}^{\top}}{x^1}+\nonumber  \\
& & \left[ 
\begin{array}{ll}
\displaystyle \min_{{x^2}(\xi^{1})}& \left[c^{{2}^{\top}}\left(\xi^{1}\right) {x^2}(\xi^{1})+ \dots + \!
										 \min_{x^H\!\left(\xxi^{H-1}\!\right)} {c^H}^{\top}\!\!\left({\xxi}^{H-1}\right) x^H\!\left(\xxi^{H-1}\!\right)\!\right]  
									\\
\text{s.t.}&   T^1(\xi^{1}){x^1} + W^2(\xi^{1}) {x^2}(\xi^{1})= h^2(\xi^{1}) \nonumber\\
&   \qquad \qquad \vdots \nonumber\\
&    T^{H-1}(\xxi^{H-1}) x^{H-1}(\xxi^{H-2}) + W^H(\xxi^{H-1}) x^H(\xxi^{H-1})= h^H(\xxi^{H-1})\nonumber\\
&   x^{t}(\xxi^{t-1})\geq 0\ ,\  t=2,\dots,H\nonumber
\end{array}\!\! \right]\!\le \gamma,\  \forall \xxi^{H-1}\in\Xi, 
\nonumber 
\end{eqnarray}
or, equivalently, as
\begin{eqnarray}
\label{ROH_X}
 \textrm{RO}_H &=&\!  \min_{{x^1},\gamma}\gamma  \\
& &\textrm{s.t. }  A {x^1}=h^1,\quad {x^1} \geq 0 \nonumber\\
& &	\qquad ({x^1},\gamma) \in \mathcal{X}_{\mathrm{RO}_H}(\xxi^{H-1}),\  \forall \xxi^{H-1} \in \Xi, \nonumber
\end{eqnarray}
where the set $\mathcal{X}_{\mathrm{RO}_H}(\xxi^{H-1})$ is defined as
\vspace{0.2cm}

$\mathcal{X}_{\mathrm{RO}_H}(\xxi^{H-1}):=$\\
\vspace{0.2cm}
\begin{displaymath}
\left\{
\begin{array}{l} 
({x^1},\gamma)\in\mathbb{R}^{n_1+1}_+ \text{ s.t. }\nonumber\\
\left[
\begin{array}{ll}
\displaystyle \min_{{x^2}(\xi^{1}),\dots,x^H\!\left(\xxi^{H-1}\right)} & c^{{1}^{\top}}{x^1} + c^{{2}^{\top}}\left(\xi^{1}\right) {x^2}(\xi^{1})+ \dots + 
\! {c^H}^{\top}\!\!\left({\xxi}^{H-1}\right)  x^H\!\left(\xxi^{H-1}\!\right)\! 	\\
\text{s.t. } &   T^1(\xi^{1}){x^1} + W^2(\xi^{1}) {x^2}(\xi^{1})= h^2(\xi^{1}) \nonumber\\
  & \qquad \qquad \vdots \nonumber\\
  &  T^{H-1}(\xxi^{H-1}) x^{H-1}(\xxi^{H-2}) + W^H(\xxi^{H-1}) x^H(\xxi^{H-1})= h^H(\xxi^{H-1})\nonumber\\
 & x^{t}(\xxi^{t-1})\geq 0\ ,\  t=2,\dots,H
\end{array} \right]\le \gamma
\end{array}
\right\}\ .
\end{displaymath}

Similarly, problem RwC$_H$ rewrites
\begin{eqnarray}
\label{RwC_XH}
 \textrm{RwC}_H &=&\!  \min_{{x^1},\gamma}\gamma  \\
& &\textrm{s.t. }  A {x^1}=h^1,\quad {x^1} \geq 0 \nonumber\\
& &	\qquad ({x^1},\gamma) \in \mathcal{X}_{\mathrm{RwC}_H}(\xxi^{H-1}),\  \forall \xxi^{H-1} \in \Xi, \nonumber
\end{eqnarray}
where the set $\mathcal{X}_{\mathrm{RwC}_H}(\xxi^{H-1})$ is defined as

\vspace{0.2cm}
$\mathcal{X}_{\mathrm{RwC}_H}(\xxi^{H-1}):=$
\vspace{0.2cm}
\begin{displaymath}
\left\{ 
({x^1},\gamma)\in\mathbb{R}^{n_1+1}_+ \text{ s.t. }
\left\{\begin{array}{l}
 \exists {x^t}(\xxi^{t-1})\in\mathbb{R}^{n_t}_+,\ t=2,\dots,H\textrm{ satisfying  }\nonumber \\ 
\displaystyle  c^{{1}^{\top}}{x^1} + c^{{2}^{\top}}\left(\xi^{1}\right) {x^2}(\xi^{1})+ \dots + 
\! {c^H}^{\top}\!\!\left({\xxi}^{H-1}\right)  x^H\!\left(\xxi^{H-1}\!\right)\leq \gamma\! 	\\
 T^1(\xi^{1}){x^1} + W^2(\xi^{1}) {x^2}(\xi^{1})= h^2(\xi^{1}) \nonumber\\
  \qquad \qquad \vdots \nonumber\\
   T^{H-1}(\xxi^{H-1}) x^{H-1}(\xxi^{H-2}) + W^H(\xxi^{H-1}) x^H(\xxi^{H-1})= h^H(\xxi^{H-1})\nonumber\
\end{array} 
\right.
\right\} \ .
\end{displaymath}
So, we just need to prove that $\mathcal{X}_{\textrm{RO}_H}(\xxi^{H-1})\equiv \mathcal{X}_{\textrm{RwC}_H}(\xxi^{H-1})$ for the minimum value of $\gamma$.
\noindent\begin{itemize}
\item[$\bullet$] We prove that if $({x^1},\gamma)\in \mathcal{X}_{\mathrm{RO_H}}$, then $({x^1},\gamma)\in \mathcal{X}_{\mathrm{RwC_H}}$. 
If $({x^1},\gamma)\in \mathcal{X}_{\textrm{RO}_H}$, then $\exists x^{t}(\xxi^{t-1})\in \mathbb{R}_{+}^{n_t}$, $t=2,\dots,H$ such that $T^{t-1}(\xxi^{t-1}){x^{t-1}}(\xxi^{t-1}) + W^t(\xxi^{t-1}){x^t}(\xxi^{t-1})= h^t(\xxi^{t-1})$, $t=2,\dots,H$ are satisfied and   
$$\min_{ {x^2}(\xi^{1}),\dots,{x^H}(\xxi^{H-1})}c^{{1}^{\top}}{x^1}+c^{{2}^{\top}}\left(\xi^{1}\right) {x^2}(\xi^{1}) + \dots + 
\! {c^H}^{\top}\!\!\left({\xxi}^{H-1}\right)  x^H\!\left(\xxi^{H-1}\!\right) \leq \gamma,$$ for the minimum value of $\gamma$.  Consequently $({x^1},\gamma)\in \mathcal{X}_{\textrm{RwC}_H}$.
\item[$\bullet$] Conversely if $({x^1},\gamma)\in \mathcal{X}_{\mathrm{RwC}_H}$, then we need to prove that $({x^1},\gamma)\in \mathcal{X}_{\mathrm{RO}_H}$.
If $({x^1},\gamma)\in \mathcal{X}_{\textrm{RwC}_H}$ then $\exists x^{t}(\xxi^{t-1})\in \mathbb{R}_{+}^{n_t}$, $t=2,\dots,H$ such that 
$T^{t-1}(\xxi^{t-1}){x^{t-1}}(\xxi^{t-1}) + W^t(\xxi^{t-1}){x^t}(\xxi^{t-1})= h^t(\xxi^{t-1})$, $t=2,\dots,H$ are satisfied and   
$c^{{1}^{\top}}{x^1}+c^{{2}^{\top}}\left(\xi^{1}\right) {x^2}(\xi^{1}) + \dots + 
\! {c^H}^{\top}\!\!\left({\xxi}^{H-1}\right)  x^H\!\left(\xxi^{H-1}\!\right) \leq \gamma$ for the minimum value of $\gamma$.

 This implies that ${x^t}(\xxi^{t-1})$, $t=2,\dots,H$ is the minimum of   $c^{{1}^{\top}}{x^1}+c^{{2}^{\top}}\left(\xi^{1}\right) {x^2}(\xi^{1}) + \dots + 
\! {c^H}^{\top}\!\!\left({\xxi}^{H-1}\right)  x^H\!\left(\xxi^{H-1}\!\right)$. By contradiction if ${x^t}(\xxi^{t-1})$, $t=2,\dots,H$ were not be the minimum then $\gamma$ would not be at the minimum of problem $\textrm{RwC}_H$. 
\qed
\end{itemize}

\section{Proof of Theorem 4}
We first prove the convexity of the set $\mathcal{X}_{\textrm{RwC}_H}(\xxi^{H-1})$ defined above.
Given $\hat{x}^1,\tilde{x}^1\in\mathcal{X}_{\textrm{RwC}_H}(\xxi^{H-1})$, then there exist $\hat{x}^t(\xxi^{t-1}),\tilde{x}^t(\xxi^{t-1})$, $t=2,\dots,H$ such that
\begin{displaymath}
 \left\{
\begin{array}{l}
 T^{t-1}(\xxi^{t-1})\hat{x}^{t-1}(\xxi^{t-2}) + W^t(\xxi^{t-1}) \hat{x}^{t}(\xxi^{t-1})= h^t(\xxi^{t-1}),\   t=2,\dots,H \\
T^{t-1}(\xxi^{t-1})\tilde{x}^{t-1}(\xxi^{t-2}) + W^t(\xxi^{t-1}) \tilde{x}^{t}(\xxi^{t-1})= h^t(\xxi^{t-1}),\   t=2,\dots,H \\
 {c^1}^{\top} \hat{x}^1+{c^{2}}^\top\!\!\!\left(\xi^{1}\right) \hat{x}^{2}\!\left(\xi^{1}\right)+\cdots+
{c^{H}}^\top\!\!\!\left(\xxi^{H-1}\right) \hat{x}^{H}\!\left(\xxi^{H-1}\right) \le \gamma \\
 {c^1}^{\top} \tilde{x}^1+{c^{2}}^\top\!\!\!\left(\xi^{1}\right) \tilde{x}^{2}\!\left(\xi^{1}\right)+\cdots+
{c^{H}}^\top\!\!\!\left(\xxi^{H-1}\right) \tilde{x}^{H}\!\left(\xxi^{H-1}\right) \le \gamma. 
\label{convexitymulti}
\end{array}
\right.
\end{displaymath}
Consider now $x^{1\lambda}:=\lambda \hat{x}^1 + (1 - \lambda)\tilde{x}^1,$ with $\lambda\in[0,1]$, and let $x^{t\lambda}:=\lambda \hat{x}^t(\xxi^{t-1}) + (1 -\lambda)\tilde{x}^t(\xxi^{t-1})$, $t=2,\dots,H$. 
We have
\begin{eqnarray}
 T^{t-1}(\xxi^{t-1})x^{t-1\lambda} + W^t(\xxi^{t-1}) x^{t\lambda} & = & T^{t-1}(\xxi^{t-1})\left(\lambda \hat{x}^{t-1}(\xxi^{t-2}) + (1 -\lambda)\tilde{x}^{t-1}(\xxi^{t-2})\right) \nonumber\\
&  & +  W^t(\xxi^{t-1}) \left(\lambda \hat{x}^t(\xxi^{t-1}) + (1 -\lambda)\tilde{x}^t(\xxi^{t-1})\right) \nonumber\\
& = & \lambda \left( T^{t-1}(\xxi^{t-1})\hat{x}^{t-1}(\xxi^{t-2}) + W^t(\xxi^{t-1}) \hat{x}^{t}(\xxi^{t-1})\right)\nonumber\\
& & +  (1 -\lambda)\left(T^{t-1}(\xxi^{t-1})\tilde{x}^{t-1}(\xxi^{t-2}) + W^t(\xxi^{t-1}) \tilde{x}^{t}(\xxi^{t-1})\right) \nonumber\\
& = &  \lambda  h^t(\xxi^{t-1}) + (1 -\lambda)h^t(\xxi^{t-1}) \nonumber\\
& = & h^t(\xxi^{t-1}),\ t=2,\dots,H,\nonumber
\end{eqnarray}
and 
\begin{eqnarray}
 & & {c^1}^{\top} x^{1\lambda}+{c^{2}}^\top\!\!\!\left(\xi^{1}\right) x^{2\lambda}+\cdots+
{c^{H}}^\top\!\!\!\left(\xxi^{H-1}\right) x^{H\lambda} \nonumber \\
&  & = {c^1}^{\top} \left(\lambda \hat{x}^1 + (1 - \lambda)\tilde{x}^1\right)\nonumber\\
& & 	\quad + {c^{2}}^\top\!\!\!\left(\xi^{1}\right)\left(\lambda \hat{x}^2(\xi^{1}) + (1 -\lambda)\tilde{x}^2(\xi^{1})\right)+\cdots\nonumber\\
& & \quad + {c^{H}}^\top\!\!\!\left(\xxi^{H-1}\right) \left(\lambda \hat{x}^H(\xxi^{H-1}) + (1 -\lambda)\tilde{x}^H(\xxi^{H-1})\right)\nonumber\\
&  & = \lambda\left({c^1}^{\top} \hat{x}^1+{c^{2}}^\top\!\!\!\left(\xi^{1}\right) \hat{x}^{2}\!\left(\xi^{1}\right)+\cdots+
{c^{H}}^\top\!\!\!\left(\xxi^{H-1}\right) \hat{x}^{H}\!\left(\xxi^{H-1}\right)\right) \nonumber\\
& & \quad + (1 -\lambda)\left({c^1}^{\top} \tilde{x}^1+{c^{2}}^\top\!\!\!\left(\xi^{1}\right) \tilde{x}^{2}\!\left(\xi^{1}\right)+\cdots+
{c^{H}}^\top\!\!\!\left(\xxi^{H-1}\right) \tilde{x}^{H}\!\left(\xxi^{H-1}\right) \right)\nonumber\\
&  &  \leq \lambda \gamma + (1 - \lambda)\gamma =\gamma, \nonumber 
\end{eqnarray}
which proves the convexity of the set $\mathcal{X}_{\textrm{RwC}_H}(\xxi^{H-1})$.
From Theorem \ref{theomultistageequiv} we observe that the condition $({x^1},\gamma)\in\mathcal{X}_{\textrm{RwC}_H}(\xxi^{H-1})$ is equivalent to $({x^1},\gamma)\in\mathcal{X}_{\textrm{RO}_H}(\xxi^{H-1})$ and that the problem $\textrm{RwC}_H$ is equivalent to $\textrm{RO}_H$ given by
\begin{eqnarray}
\label{ro_multidris}
\min_{ {x^1},\gamma} &&  \gamma  \\
& \textrm{s.t. }&			 A {x^1}=h^1,\quad {x^1}\geq 0\nonumber\\
& & c^{{1}^{\top}}{x^1} + \min_{{x^2}(\xi^1)} \left[c^{{2}^{\top}}\left(\xi^{1}\right){x^2}(\xi^1)+\dots+ \min_{x^{H}(\xxi^{H-1})} c^{{H}^{\top}}\left(\xxi^{H-1}\right){x^H}(\xxi^{H-1})\right] \leq  \gamma \nonumber \\
& & T^1({\xi}^1)x^1 + W^2(\xi^1) x^2({\xi}^1)= h^2({\xi}^1),\ \forall\xi^1\in\Xi^1\nonumber\\
& &  \qquad \qquad \vdots \nonumber\\
& &   T^{H-1}(\xxi^{H-1}) x^{H-1}(\xxi^{H-2}) + W^H(\xxi^{H-1}) x^H(\xxi^{H-1})= h^H(\xxi^{H-1}), \forall \xxi^{H-1}\in\Xi\nonumber\\
& &  x^{t}(\xxi^{t-1})\geq 0\ ,\  t=2,\dots,H\nonumber, \quad	\forall \xxi^{t-1}\in {\sf X}_{\tau=1}^{t-1} \Xi^{\tau}\nonumber.
\end{eqnarray}
For the convexity of $\mathcal{X}_{\textrm{RwC}_H}(\xxi^{H-1})$ it follows that the following functions computed at the optimal sequence of $(x^2(\xi^1),\dots,x^{H}(\xxi^{H-1}))$, say $({x^2}^{\ast}(\xi^1),\dots,{x^{H}}^{\ast}(\xxi^{H-1}))$,  $$ T^{t-1}(\xxi^{t-1}){{x^{t-1}}^{\ast}}(\xxi^{t-2}) + W^t(\xxi^{t-1}) {x^t}^{\ast}(\xxi^{t-1})= h^t(\xxi^{t}),$$ and $$c^{{1}^{\top}}{x^1} +  \left[c^{{2}^{\top}}\left(\xi^{1}\right){x^2}^{\ast}(\xi^1)+\dots+  c^{{H}^{\top}}\left(\xxi^{H-1}\right){x^H}^{\ast}(\xxi^{H-1})\right] \leq  \gamma,$$ are convex in $ {x^1}$ for given $\xxi^{H-1}$. Hence, the problem (\ref{ro_multidris}) is a robust convex optimization problem. 
Then, we construct its scenario counterpart
\begin{eqnarray}
\label{ro_munti5}
\min_{ {x^1},\gamma}&&  \gamma  \\
& \textrm{s.t. }&			  A {x^1}=h^1,\quad {x^1}\geq 0\nonumber\\
& & c^{{1}^{\top}}{x^1} + \min_{{x^2_i}} \left[ c^{{2}^{\top}}\left({\xi^{1}}^{(i)}\right){x^2_i}+\dots+ \min_{x^{H}_i} c^{{H}^{\top}}\left({\xxi^{H-1}}^{(i)}\right){x^H_i} \right]\leq  \gamma\nonumber \\
& &  T^1({\xi^{1}}^{(i)}){x^1} + W^2({\xi^{1}}^{(i)}) {x^2_i}= h^2({\xi^{1}}^{(i)}),\quad \quad {x^2_i}\geq 0,\quad	i=1,\dots,N \nonumber\\
& &  \qquad \qquad \vdots \nonumber\\
& &   T^{H-1}({\xxi^{H-1}}^{(i)}) x^{H-1}_i + W^H({\xxi^{H-1}}^{(i)}) x^H_i= h^H({\xxi^{H-1}}^{(i)}),\  i=1,\dots,N \nonumber\\
& &  x^{1}\geq 0\ ,\quad x^{t}_i\geq 0\ ,\  t=2,\dots,H,\   i=1,\dots,N\nonumber,
\end{eqnarray}
where the subscript $i$ for the variables ${x^t_i}$ highlights that the different minimization problems are independent. Finally,  we note that the  problem (\ref{ro_munti5}) corresponds to the  problem $\textrm{SwC}^N_H$ and the thesis follows from \cite{doi:10.1137/07069821X}. \qed


%
%

\begin{acknowledgements}
The authors would like to thank Daniel Kuhn and Phebe Vayanos  for helpful discussions on the inventory management with cumulative orders problem.
\end{acknowledgements}

\bibliographystyle{spmpsci} 

\end{document}